\documentclass[11pt,reqno]{amsart}
\usepackage{amssymb, amsmath, amsthm, amsfonts, amscd}
\usepackage{xcolor}
\usepackage[colorlinks=true,linkcolor=mygreen,filecolor=mybrown,citecolor=mygreen]{hyperref}
\definecolor{mygreen}{RGB}{0,128,0}
\definecolor{mybrown}{RGB}{165,42,42}
\usepackage{graphicx}
\usepackage[utf8]{inputenc}
\usepackage[english]{babel}
\usepackage{enumerate}
\usepackage{url}
\usepackage[square,sort,comma,numbers]{natbib}
\usepackage{tikz}
\usetikzlibrary{decorations.markings, arrows.meta}
\usepackage{multicol}
\usepackage{tikz}
\usepackage{geometry}
\calclayout
\usepackage{multicol}

\theoremstyle{definition}
\newtheorem{definition}{Definition}[section]
\newtheorem{theorem}{Theorem}
\newcommand{\seqnum}[1]{\text{Sequence #1}}
\usepackage[T1]{fontenc}

\usepackage{soul}
\usepackage{xcolor}
\address{Department of Mathematics, Faculty of Science, University of Bahrain, Kingdom of Bahrain}
\email{mimohamed@uob.edu.bh}

\begin{document}
	
	\title{On the Emergence of the Quanta Prime Sequence}
	\author{Moustafa Ibrahim}
	
	\maketitle
	
\begin{abstract}This paper presents the Quanta Prime Sequence (QPS) and its foundational theorem, showcasing a unique class of polynomials with substantial implications. The study uncovers profound connections between Quanta Prime numbers and essential sequences in number theory and cryptography. The investigation highlights the sequence's contribution to the emergence of new primes and its embodiment of core mathematical constructs, including Mersenne numbers, Fermat numbers, Lucas numbers, Fibonacci numbers, the Chebyshev sequence, and the Dickson sequence. The comprehensive analysis emphasizes the sequence's intrinsic relevance to the Lucas-Lehmer primality test. This research positions the Quanta Prime sequence as a pivotal tool in cryptographic applications, offering novel representations of critical mathematical structures. Additionally, a new result linking the Quanta Prime sequence to the Harmonic series is introduced, hinting at potential progress in understanding the Riemann Hypothesis.
	
\end{abstract}

	\section{Introduction}

This paper introduces the "Quanta Prime Sequence" (QPS), a mathematical construct that unifies several fundamental sequences and concepts within number theory. Defined systematically, the Quanta Prime Sequence reveals natural connections to well-known constructs such as Mersenne primes, Fibonacci numbers, Lucas numbers, Fermat numbers, and perfect numbers. This sequence not only enriches the field of number theory but also uncovers significant ties to essential mathematical tools and theories, including the Lucas-Lehmer Primality Test and harmonic numbers. These connections extend its relevance into broader mathematical and physical contexts, notably cryptography and computational mathematics.

An important attribute of the Quanta Prime Sequence is its alignment with Dickson and Chebyshev polynomials—mathematical structures pivotal to various fields such as CDMA (Code Division Multiple Access), Dickson cryptography, and permutation polynomials (see \cite{Levine}). These relationships suggest potential applications of the Quanta Prime Sequence in secure communications and signal processing, highlighting its practical significance beyond theoretical exploration.

Furthermore, the Quanta Prime Sequence's connection to the harmonic series situates it within the realm of the Riemann Hypothesis, suggesting that its properties may yield new insights into the distribution of prime numbers and the intricate structure underlying number theory. By examining these connections, this study aims to shed light on how the Quanta Prime Sequence can contribute to understanding prime number theory and related unsolved problems.

	\section{\textbf{Motivation}}
When delving into the Eight Levels Theorem, as discussed in \cite{2}, the necessity for generalization became apparent, leading to the exploration presented in \cite{3}. This exploration resulted in the discovery of a family of sequences with unique properties, denoted as $\Psi$ \cite{3}. Surprisingly, further investigation unveiled yet another distinctive and original family of sequences, termed \lq $\Omega$\rq, which is the focal point of the current paper. These sequences exhibit even more peculiar characteristics and are enriched with distinct properties related to fundamental mathematical structures and well-known arithmetic sequences, along with the emergence of a new prime. This offers diverse perspectives on essential sequences in the realms of number theory and mathematics as a whole. The present research paper sheds light on the origin of this unique and remarkable family, which we shall refer to as the Quanta Prime Sequence. In this paper, we introduce and study this new concept, along with its foundational properties. One of the key contributions of this work is the proof of the following theorem, which captures a fundamental expansion associated with the Quanta Prime Sequence.

\begin{theorem}{(The First Fundamental Theorem of the Quanta Prime sequence)}

		\label{F11} 
		For any numbers $a,b,\alpha,\beta, n$, $\beta a - \alpha b \neq 0$, we get the following expansion  		
		
		\begin{equation}
			\label{F11000} 
			\begin{aligned}
				&\quad \quad \Psi\left( \begin{array}{cc|r} a & b & n \\ \alpha & \beta & k \end{array} \right)   \\
				&=  \sum_{r=0}^{\lfloor{\frac{n}{2}}\rfloor - k} (-1)^{r+k} \: \frac{\: (n-r-k-1)! \: \:n \: \:}{(n-2r)! \: r!}
				\left(\begin{array}{c} \lfloor{\frac{n}{2}}\rfloor - r \\ k \end{array}\right)   \:	\Omega_r\big(k|\:\alpha, \beta \: | n \big) \:\:
				a^{r} \: (2a-b)^{\lfloor{\frac{n}{2}}\rfloor -k -r},
			\end{aligned}
		\end{equation}
		where the coefficients 
		\begin{equation}
			\label{F22} 
			\begin{aligned}
				(-1)^{r+k} \: \frac{\: (n-r-k-1)! \: \:n \: \:}{(n-2r)! \: r!}
				\left(\begin{array}{c} \lfloor{\frac{n}{2}}\rfloor - r \\ k \end{array}\right)   \:	\Omega_r\big(k|\:\alpha, \beta \: | n \big) 
			\end{aligned}
		\end{equation}
		are integers. 	
	\end{theorem}	
	\section{Summarized Results from \cite{3} Essential to Our Study}
	
	For a natural number $n$, we define $\delta(n) = n \pmod{2}$. For an arbitrary real number $x$, $\lfloor \frac{x}{2} \rfloor$ represents the greatest integer less than or equal to $\frac{x}{2}$. In the study conducted by \cite{3}, we summarize the main findings that are essential for the current paper, as several key results from that research will be used here.

	\begin{theorem}{(The $\Psi-$representation for sums of powers)}
		\label{WW3}
		For any natural number $n$, the $\Psi-$polynomial satisfy the following identity
		\begin{equation}
			\label{WW4}
			\begin{aligned}
				\Psi(xy,-x^2-y^2,n) &= \frac{x^n+y^n}{(x+y)^{\delta(n)}}.   \\
			\end{aligned}
		\end{equation}
	\end{theorem}

		\begin{definition}
		For any given variables $a,b$, $(a,b) \neq (0,0)$, and for any natural number $n$, we define the sequence $\Psi(a,b,n)= \Psi(n),$ by the following recurrence relation
		\begin{equation}
			\begin{aligned}
				\label{def0}
				\Psi(0)=2, \Psi(1)=1,\Psi(n+1)=(2a-b)^{\delta(n)}\Psi(n) - a \Psi(n-1).
			\end{aligned}
		\end{equation}
	\end{definition}

		\begin{theorem}
		\label{exp1}
		For any natural number $n$, and any real numbers $a,b, \alpha, \beta$, $ \beta a - \alpha b \neq 0 $, there exist unique polynomials in $a,b, \alpha, \beta$ with integer coefficients, that we call
		$ \Psi\left( \begin{array}{cc|r} a & b & n \\ \alpha & \beta & r \end{array} \right)$, that depend only on $a,b, \alpha, \beta, n,$ and $r$, and satisfy the following polynomial identity
		\begin{equation}
			\label{ex00} 
			\begin{aligned}
				(\beta a - \alpha b)^{\lfloor{\frac{n}{2}}\rfloor} \frac{x^n+y^n}{(x+y)^{\delta(n)}}  = \sum_{r=0}^{\lfloor{\frac{n}{2}}\rfloor}
				\Psi\left(\begin{array}{cc|r} a & b & n \\ \alpha & \beta & r \end{array}\right)
				(\alpha x^2 + \beta xy + \alpha y^2)^{\lfloor{\frac{n}{2}}\rfloor -r} (ax^2+bxy+ay^{2})^{r}. 
			\end{aligned}
		\end{equation}
		Moreover  
		\begin{equation}
			\label{ex000} 
			\begin{aligned}
				\Psi\left(\begin{array}{cc|r} a & b & n \\ \alpha & \beta & 0 \end{array} \right) = \Psi(a,b,n), 
			\end{aligned}
		\end{equation}
		and
		\begin{equation}
			\label{ex111}
			\begin{aligned}			
				\quad \quad \quad\Psi\left(\begin{array}{cc|c} a & b & n \\ \alpha & \beta & \lfloor{\frac{n}{2}}\rfloor \end{array} \right) = (-1)^{\lfloor{\frac{n}{2}}\rfloor} \: \Psi(\alpha,\beta,n). 
			\end{aligned}
		\end{equation}
	\end{theorem}
\begin{theorem}{(The Fundamental Theorem of the $\Psi$-Sequence)}\\
	\label{IAexp2} Let $a, b, \alpha, \beta, u, v$ be any real numbers such that $\beta a - \alpha b \neq 0$, and let $n$ be any natural number. Then we have:
	\begin{align}
		\frac{1}{\left(\left\lfloor \frac{n}{2} \right\rfloor\right)!} \left(\alpha \frac{\partial}{\partial a} + \beta \frac{\partial}{\partial b}\right)^{\left\lfloor \frac{n}{2} \right\rfloor} \Psi(a, b, n) = \Psi(\alpha, \beta, n),
	\end{align}
	where the sequence $\Psi(u, v, m) := \Psi(m)$ is defined by the recurrence relation:
	\begin{equation}
		\label{Idef0}
		\Psi(m+1) = (2u - v)^{\delta(m)} \Psi(m) - u \Psi(m-1), \quad \Psi(0) = 2, \; \Psi(1) = 1.
	\end{equation}
\end{theorem}

Furthermore, from reference \cite{3}, we have the following explicit formulas:
\begin{equation}
	\label{comp1}
	\Psi(a, b, n) = \frac{(2a - b)^{\left\lfloor \frac{n}{2} \right\rfloor}}{2^n} \left\{ \left(1 + \sqrt{\frac{b + 2a}{b - 2a}}\right)^n + \left(1 - \sqrt{\frac{b + 2a}{b - 2a}}\right)^n \right\},
\end{equation}
and
\begin{equation}
	\label{00}
	x^n + y^n = \sum_{i=0}^{\left\lfloor \frac{n}{2} \right\rfloor} (-1)^i \frac{n}{n - i} \binom{n - i}{i} (xy)^i (x + y)^{n - 2i}.
\end{equation}

\begin{theorem}
	\label{comp2}
	For any natural number $n$, the following formula is true
	\begin{equation}
		\label{comp3}
		\Psi(a,b,n) =\sum_{i=0}^{\left\lfloor \frac{n}{2} \right\rfloor}\frac{n}{n-i} \binom{n-i}{i} (-a)^i (2a-b)^{\left\lfloor \frac{n}{2} \right\rfloor - i}.
	\end{equation}
\end{theorem}

\begin{theorem}{(The product of $\Psi$-sequences)} \\
	\label{productformula}
	For any natural numbers $n,m$, the $\Psi-$polynomial satisfy the following identity
	\begin{equation}
		\label{WW8}
		\begin{aligned}
			(2a-b)^{\delta(n)\delta(m)} \: \Psi(a,b,n) \Psi(a,b,m) = \Psi(a,b,n+m) + a^{\min\{n, m\}} \Psi(a,b,n-m).   \\
		\end{aligned}
	\end{equation}
\end{theorem}

\begin{theorem}
	\label{Aexp1} For $\beta a - \alpha b \neq 0$,
	the polynomials $\Psi_r(n):=\Psi\left(\begin{array}{cc|c}
		a & b & n \\ \alpha & \beta & r \end{array}\right)$ satisfy
	\begin{align}
		\begin{aligned}
			\label{diff1}
			\big(\alpha \frac{{\partial} }{\partial a} +  \beta \frac{{\partial} }{\partial b}\big)\Psi_r(n) &= - (r+1)\Psi_{r+1}(n), \\
			\big(a \frac{{\partial} }{\partial \alpha} +  b \frac{{\partial} }{\partial \beta} \big) \Psi_r(n) &= - \big(\lfloor{\frac{n}{2}}\rfloor - r + 1 \big)\Psi_{r-1}(n).
		\end{aligned}
	\end{align}
\end{theorem}

\begin{theorem}
	\label{A1} For $\beta a - \alpha b \neq 0$,
	the polynomials $\Psi\left(\begin{array}{cc|c}
		a & b & n \\ \alpha & \beta & r \end{array}\right)$ satisfy
	\begin{align}
		\begin{aligned}
			\label{diff3}
			\Psi\left(\begin{array}{cc|c}
				a & b & n \\ \alpha & \beta & r \end{array} \right) &=  \frac{(-1)^r}{r!}
			\Big(\alpha \frac{{\partial} }{\partial a} + \beta \frac{{\partial}}{\partial b}\Big)^{r} \Psi(a,b,n),  \quad \qquad
		\end{aligned}
	\end{align}
	and
	\begin{align}
		\begin{aligned}
			\label{diff5}
			\Psi\left( \begin{array}{cc|c}
				a & b & n \\ \alpha & \beta & r \end{array} \right) =  \frac{(-1)^r}{(\lfloor{\frac{n}{2}}\rfloor -r)!} \Big(a \frac{{\partial} }{\partial \alpha} + b \frac{{\partial}}{\partial \beta}\Big)^{\lfloor{\frac{n}{2}}\rfloor-r} \Psi(\alpha,\beta,n).  
		\end{aligned}
	\end{align}
\end{theorem}

Due to the work of Euclid and Euler, it is well-known that an even integer $n$ is perfect if and only if $n = 2^{p-1}(2^p-1)$, where $2^{p}-1$ is prime. A prime of the form $2^p-1$ is called a \textit{Mersenne} prime. Recently, the paper \cite{3} mentioned the following result without proof:

\begin{theorem}{}
	\label{Theorem of ABC1}
	For any given prime $p\geq 5$, $n:=2^{p-1}$, we associate the double-indexed polynomial sequences $A_r(k)$, $B_r(k)$, which are defined by
	\begin{equation}
		\label{ABC2} 
		\begin{aligned}		
			A_r(k) &=  (p-r-k) \: A_r(k-1) + \: 4 \: (p-2r) \: A_{r+1}(k-1), \quad &A_r(0) = 1  \quad   \text{for all} \:\: r,\\
			B_r(k) &=  -2 \: (n-r-k) \: B_r(k-1) - \:2 \: (n-2r-1) \: B_{r+1}(k-1),  \quad &B_r(0) = 1  \quad   \text{for all} \:\: r.
		\end{aligned}
	\end{equation}
	Then both of the ratios   
	\begin{equation}
		\label{ABC3} 
		\begin{aligned}
			\frac{A_0( \lfloor{\frac{p}{2}}\rfloor   )}{(p-1)(p-2) \cdots (p - \lfloor{\frac{p}{2}}\rfloor )} \quad	, \quad 	\frac{B_0( \lfloor{\frac{n}{2}}\rfloor   )}{(n-1)(n-2) \cdots (n - \lfloor{\frac{n}{2}}\rfloor )} 	
		\end{aligned}
	\end{equation}
	are integers. Moreover the number $2^p -1$ is prime \bf{if and only if} 
	\begin{equation}
		\label{ABC4} 
		\begin{aligned}
			\frac{A_0( \lfloor{\frac{p}{2}}\rfloor   )}{(p-1)(p-2) \cdots (p - \lfloor{\frac{p}{2}}\rfloor )}  \quad  \Big\vert \quad \frac{B_0( \lfloor{\frac{n}{2}}\rfloor   )}{(n-1)(n-2) \cdots (n - \lfloor{\frac{n}{2}}\rfloor )}.
		\end{aligned}
	\end{equation}
	
\end{theorem}

\section{\textbf{Objective of the Current Paper}}
This paper aims to provide a comprehensive introduction to the Quanta Prime sequence. Additionally, it endeavors to unveil the fundamental theorems governing the Quanta Prime sequence and elucidate its connections with various mathematical topics. Furthermore, the paper presents an inaugural and comprehensive proof for Theorem \eqref{Theorem of ABC1}, showcasing its originality and contributing another valuable generalization.

\subsection*{Key Findings}
The key contributions of this paper include:

\begin{enumerate}
	
	\item Introducing the Quanta Prime Sequence.
	
	\item Presentation of the fundamental theorem governing the Quanta Prime Sequence, detailing its significance and implications for mathematical research.
	
	\item Elucidation of connections with diverse mathematical domains, exploring how the Quanta Prime Sequence interrelates with various topics such as Mersenne primes, Fibonacci numbers, and Lucas numbers.
	
	\item Comprehensive proof of Theorem \eqref{Theorem of ABC1}, providing a thorough and rigorous demonstration of its validity and applications.
	
	\item Investigation into the emergence of new primes $p_k$, examining the role of the Quanta Prime Sequence in understanding the progression and distribution of prime numbers.
	
	\item Highlighting the connections between the Quanta Prime Sequence and Fibonacci numbers, analyzing the interplay and shared properties between these fundamental sequences.
	
	\item Highlighting, motivating, and predicting future research connections among the Quanta Prime sequence, the Harmonic series, and the Riemann Hypothesis.

\end{enumerate}

\section{\textbf{Key Discoveries Presented in this Paper}}This section provides a concise overview of the significant advancements and newly introduced theorems, along with their complete proofs, as presented in the current paper.

		\begin{theorem}{ (The Second Fundamental Theorem of the Quanta Prime Sequence)}\\
		\label{Ispace3} 
		For any natural numbers $n, k$, and any nonzero point  $(\zeta, \xi)$, we get
		\begin{equation}
			\label{Ispace4} 
			\begin{aligned}
				\Psi(\zeta, \xi, n) 	    \:= \: \frac{\:	\Omega_0\big(\lfloor{\frac{n}{2}}\rfloor|\:\zeta, \xi\: | n \big)  \:}{ \; (n-1)(n-2) \cdots (n - \lfloor{\frac{n}{2}}\rfloor ) \:     },	
			\end{aligned}
		\end{equation}	
where the sequence $\Omega$ is defined by the double-indexed recurrence relation  
\begin{equation}  	
		\begin{aligned}
		\Omega_r\big(k|\zeta, \xi|n\big) &= (2\zeta-\xi)  \: \big(n-r-k\big) \: \Omega_r\big(k-1|\zeta, \xi|n\big) \\
		& \qquad \qquad \quad - 2\ \zeta  \: \big(n-2r-\delta(n-1)\big) \: \Omega_{r+1}\big(k-1|\zeta, \xi|n\big),  \\
		\Omega_r\big(0|\zeta, \xi|n\big) &= 1   \quad   \text{for all} \:\: r.
	\end{aligned}
\end{equation}
Moreover, the $\Omega$ sequence unexpectedly satisfies the  identity
\begin{equation}
	\label{IFA2} 
	\begin{aligned}
		\frac{x^n+y^n}{(x+y)^{\delta(n)}}\:= \: \frac{\Omega_0\big(\lfloor{\frac{n}{2}}\rfloor |\:xy , - x^2 - y^2 \: |n \big) }{(n-1)(n-2) \cdots (n - \lfloor{\frac{n}{2}}\rfloor )}. 	
	\end{aligned}
\end{equation}
Moreover, it is the following expansion that we observed the existence of $\Omega$ sequence 

	\begin{equation}
	\label{IF1100} 
	\begin{aligned}
		&\quad \quad \frac{(-1)^k}{k!}
		\Big(\alpha \frac{{\partial} }{\partial a} + \beta \frac{{\partial}}{\partial b}\Big)^{k} \Psi(a,b,n)  \\
		&=  \sum_{r=0}^{\lfloor{\frac{n}{2}}\rfloor - k} (-1)^{r+k} \: \frac{\: (n-r-k-1)! \: \:n \: \:}{(n-2r)! \: r!}
		\left(\begin{array}{c} \lfloor{\frac{n}{2}}\rfloor - r \\ k \end{array}\right)   \:	\Omega_r\big(k|\:\alpha, \beta \: | n \big) \:\:
		a^{r} \: (2a-b)^{\lfloor{\frac{n}{2}}\rfloor -k -r}.
	\end{aligned}
\end{equation}
We call the $\Omega-$sequence the Quanta Prime Sequence. Moreover, for any nonzero point $(\alpha,\beta)$, we get the following property 

\begin{equation}
	\label{omar} 
	\begin{aligned}
		p_{k+1}\:  \Big\vert \: \Omega_0\big(p_{k}|\:\alpha, \beta\: | 2p_{k} \big),						
	\end{aligned}
\end{equation}	
where $p_k$ denotes the $k$-th prime, $k >1$.
\end{theorem}

		\begin{theorem}{(New representation of Mersenne numbers)}
		\label{Theorem of G2fQ}
		For any given odd natural number $p$, the number $2^p -1$ can be represented by
		\begin{equation}
			\label{G2BfQ} 
			\begin{aligned}
				2^p-1 \:=\: \frac{A_0( \lfloor{\frac{p}{2}}\rfloor   )}{(p-1)(p-2) \cdots (p - \lfloor{\frac{p}{2}}\rfloor )}, 	
			\end{aligned}
		\end{equation}
		
		where the double-indexed polynomial sequence $A_r(k)$ is defined by the recurrence relation 
		\begin{equation}
			\label{G2-e-q} 
			\begin{aligned}
				A_r(k) &=  (p-r-k) \: A_r(k-1) + \: 4 \: (p-2r) \: A_{r+1}(k-1), \quad &A_r(0) = 1  \quad   \text{for all} \:\: r.
			\end{aligned}
		\end{equation}
	\end{theorem}

	\begin{theorem}{}
		\label{KG11}
		For any given natural number $n$, we associate the double-indexed polynomial sequence $U_r(k)$, which is defined by
		\begin{equation}
			\label{KG22} 
			\begin{aligned}
				U_r(k) &=  \: (n-r-k) \: U_r(k-1) - \:2 \: (n-2r-\delta(n-1)) \: U_{r+1}(k-1), \\
				 \quad U_r(0) &= 1   \quad   \text{for all} \:\: r.
			\end{aligned}
		\end{equation}
		Then 	
		\begin{equation}
			\label{KG44} 
			\frac{U_0( \lfloor{\frac{n}{2}}\rfloor   )}{(n-1)(n-2) \cdots (n - \lfloor{\frac{n}{2}}\rfloor )} 	\:=
			\begin{cases}
				+2	      &   n \equiv \pm 0       \pmod{6} \\
				+1	      &   n \equiv \pm 1       \pmod{6} \\
				-1        &   n \equiv \pm 2       \pmod{6} \\
				-2	      &   n \equiv \pm 3      \pmod{6} 
			\end{cases}.    
		\end{equation}
	\end{theorem}
	
		\begin{theorem}{}
		\label{KG11Q}
		For any given natural number $n$, we associate the double-indexed polynomial sequence $V_r(k)$, which is defined by
		\begin{equation}
			\label{KG22Q} 
			\begin{aligned}
				V_r(k) &=  +2 \: (n-r-k) \: V_r(k-1) - \:2 \: (n-2r-\delta(n-1)) \: V_{r+1}(k-1), \\ \quad V_r(0) &= 1   \quad   \text{for all} \:\: r.
			\end{aligned}
		\end{equation}
		Then 	
		\begin{equation}
			\label{KG44Q} 
			\frac{V_0( \lfloor{\frac{n}{2}}\rfloor   )}{(n-1)(n-2) \cdots (n - \lfloor{\frac{n}{2}}\rfloor )} 	\:=
			\begin{cases}
				+2	      &   n \equiv \pm 0       \pmod{8} \\
				+1	      &   n \equiv \pm 1       \pmod{8} \\
				\: 0	      &   n \equiv \pm 2       \pmod{8} \\				
				-1        &   n \equiv \pm 3       \pmod{8} \\
				-2	      &   n \equiv \pm 4      \pmod{8} 
			\end{cases}.    
		\end{equation}
	\end{theorem}

	\begin{theorem}{}
		\label{PtP}
		For any given natural number $n$, we associate the double-indexed polynomial sequence $W_r(k)$, which is defined by
		\begin{equation}
			\label{GG1-1} 
			\begin{aligned}
				W_r(k) &=  3 \: (n-r-k) \: W_r(k-1) - \:2 \: (n-2r-\delta(n-1)) \: W_{r+1}(k-1), \\ \:\: W_r(0) &= 1   \   \text{for all} \:\: r.
			\end{aligned}
		\end{equation}
		Then 		
		\begin{equation}
			\label{PPP1} 
			\: 	\frac{W_0( \lfloor{\frac{n}{2}}\rfloor   )}{(n-1)(n-2) \cdots (n - \lfloor{\frac{n}{2}}\rfloor )} 	 \:=
			\begin{cases}
				+2	      &   n \equiv \pm 0       \pmod{12} \\
				+1	      &  n \equiv \pm 1 , \pm 2      \pmod{12} \\
				\:\:0        &   n \equiv \pm 3       \pmod{12} \\
				-1        &  n \equiv \pm 4, \pm 5       \pmod{12} \\
				-2	      &   n \equiv \pm 6  \pmod{12}
			\end{cases}.  
		\end{equation}
	\end{theorem}		
	
	\begin{theorem}
		\label{Theorem of G5}
		For any given natural number $n$, we associate the double-indexed polynomial sequence $T_r(k)$,  which is defined by
		\begin{equation}
			\label{G5} 
			\begin{aligned}
				T_r(k) &=  4 \: (n-r-k) \: T_r(k-1)  - 2\: (n-2r - \: \delta(n-1)) \: T_{r+1}(k-1),  \\
				T_r(0) &= 1   \quad   \text{for all} \:\: r.
			\end{aligned}
		\end{equation}
	Then
		\begin{equation}
			\label{G5-2} 
			\begin{aligned}
				 \frac{T_0( \lfloor{\frac{n}{2}}\rfloor   )}{(n-1)(n-2) \cdots (n - \lfloor{\frac{n}{2}}\rfloor )} \:=\:    2^{\delta(n-1)}. 	
			\end{aligned}
		\end{equation}
	\end{theorem}

	\begin{theorem}{(New representation for Lucas sequence)}
		\label{Theorem of G6}
		For any given natural number $n$, we associate the double-indexed polynomial sequence $H_r(k)$, which is defined by
		\begin{equation}
			\label{G6} 
			\begin{aligned}
				H_r(k) &=   \: (n-r-k) \: H_r(k-1) + 2\: (n-2r - \: \delta(n-1)) \: H_{r+1}(k-1), \\
					H_r(0) &= 1   \quad   \text{for all} \:\: r.
			\end{aligned}
		\end{equation}
		Then 
		\begin{equation}
			\label{G6-2} 
			\begin{aligned}
				L(n) \:= \: \frac{H_0(\lfloor{\frac{n}{2}}\rfloor   )}{(n-1)(n-2) \cdots (n - \lfloor{\frac{n}{2}}\rfloor )}, 	
			\end{aligned}
		\end{equation}
	
		where $L(n)$ is Lucas sequence defined by $L(m+1)=L(m) + L(m-1)$, $L(1)=1$, $L(0)=2$.  	\end{theorem}

	\begin{theorem}{(New representation for Fermat numbers)}
		\label{Theorem of G4}
		For any given natural number $n$, we associate the double-indexed polynomial sequence $F_r(k)$, which is defined by
		\begin{equation}
			\label{G4} 
			\begin{aligned}
				F_r(k) &=   \: \big(2^n-r-k\big) \: F_r(k-1) + 4 \: \big(2^n -2r -1 \big) \: F_{r+1}(k-1),  \\
				F_r(0) &= 1   \quad   \text{for all} \:\: r.
			\end{aligned}
		\end{equation}
		Then the Fermat number $F_{n} = 2^{2^n} + 1 $ can be represented by 
		\begin{equation}
			\label{G4-2} 
			\begin{aligned}
				F_{n}  \:=\: \frac{F_0\big(2^{n-1}   \big)}{\big(2^n-1\big)\big(2^n-2\big) \cdots \big(2^{n-1}\big)}. 	
			\end{aligned}
		\end{equation}
	\end{theorem}
	\begin{theorem}{(New representation for Fibonacci-Lucas oscillating sequence)}
		\label{Theorem of G7}
		For any given natural number $n$, we associate the double-indexed polynomial sequence $G_r(k)$, which is defined by
		\begin{equation}
			\label{G7} 
			\begin{aligned}
			G_r(k) &=   5 \: (n-r-k) \: G_r(k-1) - 2\: \big(n-2r-\delta(n-1)\big) \: G_{r+1}(k-1),  \\
				G_r(0) &= 1   \quad   \text{for all} \:\: r.
			\end{aligned}
		\end{equation}
		Then 
		\begin{equation}
			\label{G7-2} 
			\begin{aligned}
				\frac{G_0( \lfloor{\frac{n}{2}}\rfloor   )}{(n-1)(n-2) \cdots (n - \lfloor{\frac{n}{2}}\rfloor )} = \begin{cases}
					F(n)	      &   \text{if} \quad n \:\: odd \\
					L(n)	      &  \text{if} \quad n \: \: even  
				\end{cases},	
			\end{aligned}
		\end{equation}
		where $F(n)$ is Fibonacci sequence defined by $F(m+1)=F(m) + F(m-1), F(1)=1, F(0)=0$, and
	$L(n)$ is Lucas sequence defined by $L(m+1)=L(m) + L(m-1)$, $L(1)=1$, $L(0)=2$.	
		 
	\end{theorem}

	\section{\textbf{The Quanta Prime sequence}}

	\begin{definition}{(The  Quanta Prime sequence associated with $n$ and a point $(\zeta, \xi )$)   }
		\label{omegaDef} \\
		For any given natural number $n$, and a point $(\zeta, \xi)$, $(\zeta, \xi)\neq (0,0)$, we associate the double-indexed sequence $\Omega_r(k|\zeta, \xi|n)$, where $0 \leq r + k \leq  \lfloor{\frac{n}{2}}\rfloor $, such that 
		\begin{equation}
			\label{G0} 
			\begin{aligned}
	\Omega_r\big(k|\zeta, \xi|n\big) &= (2\zeta-\xi)  \: \big(n-r-k\big) \: \Omega_r\big(k-1|\zeta, \xi|n\big) \\
	& \quad - 2\ \zeta  \: \big(n-2r-\delta(n-1)\big) \: \Omega_{r+1}\big(k-1|\zeta, \xi|n\big),  \\
				\Omega_r\big(0|\zeta, \xi|n\big) &= 1   \quad   \text{for all} \:\: r.
			\end{aligned}
		\end{equation}
		
	We call $\Omega_r\big(k|\zeta, \xi|n\big)$ the Quanta Prime Sequence (QPS) associated with $n$ and the point $(\zeta, \xi)$ at level $k$ with variation $r$.
		
We define $\Omega_{\text{stable}}\big(k|\zeta, \xi|n\big)$ as the Quanta Prime Sequence associated with $n$ and the point $(\zeta, \xi)$ at level $k$, characterized by stability in the absence of variation ($r$). When the variation term is removed, the sequence maintains a constant and unchanging nature, signifying a reliable and consistent pattern across different levels and points. This stability in the Quanta Prime Sequence suggests a persistent and uniform structure, providing valuable insight into its behavior and properties in the absence of external variations. Researchers can rely on the stability of this sequence to observe its inherent characteristics and study its behavior with greater predictability, making it a key aspect of mathematical analysis and exploration.

	\end{definition}
	
	\subsection*{Notation}
		For a given point $(\zeta, \xi),$ we put \[   \Omega_r\big(k|\zeta, \xi|n\big) = \Omega_r(k|n). \]
		For a given $n$, and point $(\zeta, \xi),$ we put \[   \Omega_r\big(k|\zeta, \xi|n\big) = \Omega_r(k). \]
		
	\section{Illustrative Example: Calculating the Quanta Prime Sequence}
		\subsection*{The  Quanta Prime Sequence associated with $(1,-2)$ }
	In this section, we present a detailed example that demonstrates the step-by-step calculation of the Quanta Prime sequence. This example serves to clarify the sequence generation process, offering a practical and accessible understanding of its structure and progression.

\begin{figure}
	\begin{tikzpicture}
		\pgfmathsetmacro{\NumBalls}{4}
		
		\pgfmathsetmacro{\VerticalDistance}{2}
		\pgfmathsetmacro{\HorizontalSeparation}{5}
		\pgfmathsetmacro{\AdditionalSeparation}{2}
		
		\foreach \i in {0,...,\NumBalls} {
			\fill (0, -\i*\VerticalDistance) circle (0.2cm);
			\ifnum\i>0
			\draw[-{Latex[length=5mm,width=2mm]}, line width=1.5pt] (0, -\i*\VerticalDistance) -- (0, -\the\numexpr\i-1\relax*\VerticalDistance);
			\fi
		}
		
		\node[right=0.2cm, font=\small] at (0, -\NumBalls*\VerticalDistance) {$1$};
		\node[right=0.2cm, font=\small] at (0, -\numexpr\NumBalls-1\relax*\VerticalDistance) {$2(n-1)$};
		\node[right=0.2cm, font=\small] at (0, -\numexpr\NumBalls-2\relax*\VerticalDistance) {$2^2(n-1)(n-3)$};
		\node[right=0.2cm, font=\small] at (0, -\numexpr\NumBalls-3\relax*\VerticalDistance) {$2^3(n-1)(n-3)(n-5)$};
		
\node[right=0.2cm, font=\small] at (\HorizontalSeparation + \AdditionalSeparation, {(\numexpr\NumBalls-4\relax)*\VerticalDistance}) {\dots};

		\node[below right, align=center, font=\small] at (\HorizontalSeparation + \AdditionalSeparation - 10, -\NumBalls*\VerticalDistance - 0.5) {For $n$ even, the Levels of \\ The Quanta Prime sequence\\associated with $n$ and about $(1,-2)$};
		
		\foreach \i in {0,...,\NumBalls} {
			\fill (\HorizontalSeparation + \AdditionalSeparation, -\i*\VerticalDistance) circle (0.2cm);
			\ifnum\i>0
			\draw[-{Latex[length=5mm,width=2mm]}, line width=1.5pt] (\HorizontalSeparation + \AdditionalSeparation, -\i*\VerticalDistance) -- (\HorizontalSeparation + \AdditionalSeparation, -\the\numexpr\i-1\relax*\VerticalDistance);
			\fi
		}
		
		\node[right=0.2cm, font=\small] at (\HorizontalSeparation + \AdditionalSeparation, -\NumBalls*\VerticalDistance) {$1$};
		\node[right=0.2cm, font=\small] at (\HorizontalSeparation + \AdditionalSeparation, -\numexpr\NumBalls-1\relax*\VerticalDistance) {$2(n-2)$};
		\node[right=0.2cm, font=\small] at (\HorizontalSeparation + \AdditionalSeparation, -\numexpr\NumBalls-2\relax*\VerticalDistance) {$2^2(n-2)(n-4)$};
		\node[right=0.2cm, font=\small] at (\HorizontalSeparation + \AdditionalSeparation, -\numexpr\NumBalls-3\relax*\VerticalDistance) {$2^3(n-2)(n-4)(n-6)$};
		
\pgfmathsetmacro{\PositionX}{\HorizontalSeparation + \AdditionalSeparation}
\pgfmathsetmacro{\PositionY}{-\numexpr\NumBalls-4\relax * \VerticalDistance}
\node[right=0.2cm, font=\small] at (\PositionX, \PositionY) {$\dots$};

		\node[below right, align=center, font=\small] at (\HorizontalSeparation + \AdditionalSeparation - 3.2, -\NumBalls*\VerticalDistance - 0.5) {For $n$ odd, the Levels of \\ The Quanta Prime sequence\\associated with $n$ and about $(1,-2)$};
	\end{tikzpicture}
	\caption{The Stable Quanta Prime Sequence at $(1,-2)$}\label{Quanta}
\end{figure}
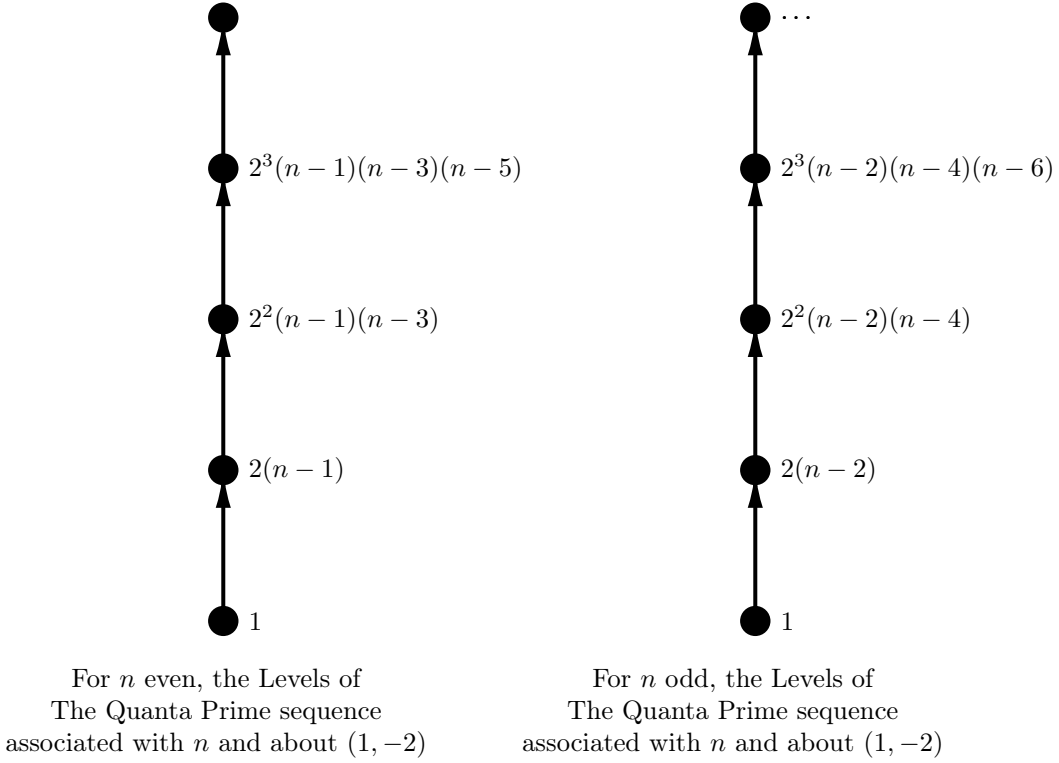

\begin{theorem}
	\label{AU5}
	The Omega sequence associated with \( n \) and the point \( (1, -2) \) is given by
	\begin{equation}
		\label{AU6} 
		\Omega_r\big(k|1, -2|n\big) = 2^k \prod\limits_{\lambda = 1}^{k}\big(n + \delta(n - 1) - 2 \lambda \big).
	\end{equation}
\end{theorem}

\begin{proof}
	We know that the Quanta Prime sequence associated with \( n \) and \( (1, -2) \) is given by the recurrence relation 
	\begin{equation}
		\label{AU1} 
		\begin{aligned}
			\Omega_r\big(k|1, -2|n\big) &= 4 \big(n - r - k\big) \Omega_r\big(k - 1|1, -2|n\big) \\
			& \quad - 2 \big(n - 2r - \delta(n - 1)\big) \Omega_{r + 1}\big(k - 1|1, -2|n\big), \\
			\Omega_r\big(0|1, -2|n\big) &= 1 \quad \text{for all} \:\: r.
		\end{aligned}
	\end{equation}
	
	Solving \eqref{AU1} iteratively, we find:
	
	\begin{equation}
		\label{AU2} 
		\begin{aligned}
			\Omega_r\big(1|1, -2|n\big) &= 4 \big(n - r - 1\big) \Omega_r\big(0|1, -2|n\big) - 2 \big(n - 2r - \delta(n - 1)\big) \Omega_{r + 1}\big(0|1, -2|n\big) \\
			&= 4 \big(n - r - 1\big) (1) - 2 \big(n - 2r - \delta(n - 1)\big) (1) \\
			&= 2 \big(n + \delta(n - 1) - 2\big).  
		\end{aligned}
	\end{equation}
	
	Continuing in this way, we find:
	
	\begin{equation}
		\label{AU3} 
		\begin{aligned}
			\Omega_r\big(2|1, -2|n\big) &= 4 \big(n - r - 2\big) \Omega_r\big(1|1, -2|n\big) - 2 \big(n - 2r - \delta(n - 1)\big) \Omega_{r + 1}\big(1|1, -2|n\big) \\
			&= 2^2 \big(n + \delta(n - 1) - 2\big) \big(n + \delta(n - 1) - 4\big). 
		\end{aligned}
	\end{equation}
	
	Repeating this process, we find:
	
	\begin{equation}
		\label{AU4} 
		\begin{aligned}
			\Omega_r\big(3|1, -2|n\big) &= 4 \big(n - r - 3\big) \Omega_r\big(2|1, -2|n\big) - 2 \big(n - 2r - \delta(n - 1)\big) \Omega_{r + 1}\big(2|1, -2|n\big) \\
			&= 2^3 \big(n + \delta(n - 1) - 2\big) \big(n + \delta(n - 1) - 4\big) \big(n + \delta(n - 1) - 6\big).
		\end{aligned}
	\end{equation}
	
	Proceeding in this way, and by induction, we can show that the general form is given by:
	
	\begin{equation}
		\label{AU7} 
		\Omega_r\big(k|1, -2|n\big) = 2^k \prod_{\lambda=1}^{k}\big(n + \delta(n - 1) - 2\lambda\big).
	\end{equation}
	
	This completes the proof of Theorem \eqref{AU5}.
\end{proof}

\begin{theorem}
	\label{AU5N}
	For a given natural number
	
	 \( n \), the following identity holds:
	\begin{equation}
		\label{AU6N}
		\Omega_0\left( \left\lfloor \frac{n}{2} \right\rfloor \mid 1, -2 \mid n \right) = 2^{\left\lfloor \frac{n}{2} \right\rfloor} \prod_{\lambda = 1}^{\left\lfloor \frac{n}{2} \right\rfloor} \left(n + \delta(n-1) - 2\lambda \right).
	\end{equation}
\end{theorem}

\begin{proof}
	Setting \( r = 0 \) and \( k = \left\lfloor \frac{n}{2} \right\rfloor \) in \eqref{AU6} directly yields the result in Theorem \ref{AU5N}.
\end{proof}

The following theorem provides an explicit formula for the sequence \( \Omega_0 \):

\begin{theorem}
	\label{AU5NM}
	\begin{equation}
		\label{AU6NM} 
		\Omega_0\left( \left\lfloor \frac{n}{2} \right\rfloor \mid 1, -2 \mid n \right) = 2^{\left\lfloor \frac{n}{2} \right\rfloor} \big(n + \delta(n - 1) - 2\big)\big(n + \delta(n - 1) - 4\big) \dotsm (1).
	\end{equation}
\end{theorem}

\begin{proof}
	Note that
	\begin{equation}
		\begin{aligned}
			n + \delta(n - 1) - 2 \left\lfloor \frac{n}{2} \right\rfloor &= n + \delta(n - 1) - 2 \cdot \frac{n - \delta(n)}{2} \\
			&= \delta(n - 1) + \delta(n) = 1.
		\end{aligned}
	\end{equation}	
	Thus, we obtain the formula presented in Theorem \ref{AU5NM}.
\end{proof}

	\section{The representation of \texorpdfstring{$\Psi$}{} in terms of \texorpdfstring{$\Omega$}{} sequence }
	
	\begin{theorem}{}
		\label{FD1} 
		For any $n,k,a,b,\alpha,\beta, \theta, n$, $\beta a - \alpha b \neq 0,$ the following expansion is true
		\begin{equation}
			\label{FD2}
			\begin{aligned}
				\Psi\left(\begin{array}{cc|c} a & b & n \\ \alpha & \beta & k \end{array} \right) = \sum_{r=0}^{\left\lfloor \frac{n}{2} \right\rfloor -k} \frac{(-1)^k}{k!}  \: \lambda_r(k|\alpha, \beta|n) a^r (2a-b)^{\left\lfloor \frac{n}{2} \right\rfloor -k- r}
			\end{aligned}
		\end{equation}
		where the numbers $\lambda_r(k|\alpha, \beta|n)$ are integers and divisible by $k!$ and satisfy the double-indexed recurrence relation 		
		
		\begin{equation}
			\label{FD3}
			\begin{aligned}
				\begin{cases}
					&        \lambda_r(k|\alpha, \beta|n)= 
					\:	(2\alpha - \beta) \:
					\Big(\left\lfloor \frac{n}{2} \right\rfloor -k- r +1 \Big)
					\: \lambda_r(k-1|\alpha, \beta|n)      \\	  
					&\qquad \qquad \qquad  \qquad \qquad \qquad  \qquad \qquad \:\: \: + \quad  \alpha \: (r+1) \: \lambda_{r+1}(k-1|\alpha, \beta|n)         	    \:     \\ 
					&\: \lambda_r(\:0\:|\:\alpha, \beta \:|\:n\:) = (-1)^r  \frac{n}{n-r} \binom{n-r}{r}   \\
				\end{cases} 
			\end{aligned}
		\end{equation}
	\end{theorem}
	
		\begin{proof}  
		From Theorem\eqref{A1}, we know that 
		\[ \Psi\left(\begin{array}{cc|c} a & b & n \\ \alpha & \beta & k \end{array} \right) =  \frac{(-1)^k}{k!} \Big(\alpha \frac{{\partial} }{\partial a} + \beta \frac{{\partial}}{\partial b}\Big)^{k} \Psi(a,b,n). \]			
		Then from equation \eqref{comp3}, it follows  	
		\begin{equation}
			\label{Q1}
			\begin{aligned}
				\Psi\left(\begin{array}{cc|c} a & b & n \\ \alpha & \beta & k \end{array} \right)  &= \frac{(-1)^k}{k!} \Big(\alpha \frac{{\partial} }{\partial a} + \beta \frac{{\partial}}{\partial b}\Big)^{k}  \sum_{r=0}^{\left\lfloor \frac{n}{2} \right\rfloor}\frac{n}{n-r} \binom{n-r}{r} (-a)^r (2a-b)^{\left\lfloor \frac{n}{2} \right\rfloor - r} \\
				&= \frac{(-1)^k}{k!} \sum_{r=0}^{\left\lfloor \frac{n}{2} \right\rfloor} (-1)^r  \frac{n}{n-r} \binom{n-r}{r}  \Big(\alpha \frac{{\partial} }{\partial a} + \beta \frac{{\partial}}{\partial b}\Big)^{k}   a^r (2a-b)^{\left\lfloor \frac{n}{2} \right\rfloor - r}. 
			\end{aligned}
		\end{equation}
		From \eqref{Q1}, and for $0 \leq k+r \leq \left\lfloor \frac{n}{2} \right\rfloor$, there exist integers $\lambda_r(k|\alpha, \beta|n)$, which  are divisible by $k!$, and depend only on the numbers $r,k,\alpha, \beta,n$ (and independent on $a,b$), such that   
		\begin{equation}
			\label{Q2}
			\begin{aligned}
				\Psi\left(\begin{array}{cc|c} a & b & n \\ \alpha & \beta & k \end{array} \right) = \frac{(-1)^k}{k!} \sum_{r=0}^{\left\lfloor \frac{n}{2} \right\rfloor -k}  \lambda_r(k|\alpha, \beta|n) a^r (2a-b)^{\left\lfloor \frac{n}{2} \right\rfloor -k- r}. \\
			\end{aligned}
		\end{equation}
		This proves \eqref{FD2}. Now, to study the coefficients $\lambda_r(k|\alpha, \beta|n)$, we need to compute the recurrence relation that arise up easily once we notice, from Theorem \eqref{Aexp1}, that the following differential property of $\Psi$ for any non-negative integer $k$ 
		
		\begin{align}
			\begin{aligned}
				\label{Q4}
				\big(\alpha \frac{{\partial} }{\partial a} +  \beta \frac{{\partial} }{\partial b}\big)	\Psi\left(\begin{array}{cc|c} a & b & n \\ \alpha & \beta & k \end{array} \right) = - (k+1)	\Psi\left(\begin{array}{cc|c} a & b & n \\ \alpha & \beta & k+1 \end{array} \right).   
			\end{aligned}
		\end{align}	
		Then from \eqref{Q2}, \eqref{Q4}, we get 	
		\begin{equation}
			\label{Q5}
			\begin{aligned}
				\big(\alpha \frac{{\partial} }{\partial a} +  \beta \frac{{\partial} }{\partial b}\big)		\frac{(-1)^k}{k!} &\sum_{r=0}^{\left\lfloor \frac{n}{2} \right\rfloor -k}  \lambda_r(k|\alpha, \beta|n) a^r (2a-b)^{\left\lfloor \frac{n}{2} \right\rfloor -k- r} \\ &= - (k+1)  \frac{(-1)^{(k+1)}}{(k+1)!} \sum_{r=0}^{\left\lfloor \frac{n}{2} \right\rfloor -k -1}  \lambda_r(k+1|\alpha, \beta|n) a^r (2a-b)^{\left\lfloor \frac{n}{2} \right\rfloor -k-1- r}.
			\end{aligned}
		\end{equation}	
		Simplifying, and noting that the coefficients $\lambda_r(k|\alpha, \beta|n)$ are independent on $a,b$, we get
		\[ \big(\alpha \frac{{\partial} }{\partial a} +  \beta \frac{{\partial} }{\partial b}\big)  \lambda_r(k|\alpha, \beta|n) =0.                  \] 
		Therefore
		\begin{equation}
			\label{Q6}
			\begin{aligned}
				\sum_{r=0}^{\left\lfloor \frac{n}{2} \right\rfloor -k}  \lambda_r(k|\alpha, \beta|n)        &  	\big(\alpha \frac{{\partial} }{\partial a} +  \beta \frac{{\partial} }{\partial b}\big)	            a^r (2a-b)^{\left\lfloor \frac{n}{2} \right\rfloor -k- r} \\ &=  \sum_{r=0}^{\left\lfloor \frac{n}{2} \right\rfloor -k -1}  \lambda_r(k+1|\alpha, \beta|n) a^r (2a-b)^{\left\lfloor \frac{n}{2} \right\rfloor -k- r -1}.
			\end{aligned}
		\end{equation}		
		Hence
		\begin{equation}
			\label{Q7}
			\begin{aligned}
				\sum_{r=0}^{\left\lfloor  \frac{n}{2} \right\rfloor -k -1}  \alpha &\: (r+1) \: \lambda_{r+1}(k|\alpha, \beta|n)         	    \:     a^r \:  (2a-b)^{\left\lfloor \frac{n}{2} \right\rfloor -k- r-1} \\ 
				+	\sum_{r=0}^{\left\lfloor \frac{n}{2} \right\rfloor -k -1} & (2\alpha - \beta)
				\Big(\left\lfloor \frac{n}{2} \right\rfloor -k- r \Big)
				\: \lambda_r(k|\alpha, \beta|n)      	    \:     a^r \: (2a-b)^{\left\lfloor \frac{n}{2} \right\rfloor -k- r-1} \\
				=  \sum_{r=0}^{\left\lfloor \frac{n}{2} \right\rfloor -k -1} & \lambda_r(k+1|\alpha, \beta|n) \: a^r \: (2a-b)^{\left\lfloor \frac{n}{2} \right\rfloor -k- r-1}.
			\end{aligned}
		\end{equation}	
		From \eqref{Q7}, comparing the coefficients, and noting that $a, 2a-b$ are algebraically independent, we immediately get 
	
		\begin{equation}
			\label{Q9}
			\begin{aligned}
				\alpha \: (r+1)  \: \lambda_{r+1}(k|\alpha, \beta|n)  
				+ \:	(2\alpha - \beta) \:
				\Big(\left\lfloor \frac{n}{2} \right\rfloor -k- r \Big)&
				\:  \lambda_r(k|\alpha, \beta|n)       	  \\
				&\qquad  =   \lambda_r(k+1|\alpha, \beta|n).
			\end{aligned}
		\end{equation}
		
		Now the initial value for $\lambda_r(\:k\:|\alpha, \beta \:|\:n\:)$, that corresponds to $k=0$, is given by
		
		\begin{equation}
			\label{FD2Q}
			\begin{aligned}
				\Psi\left(\begin{array}{cc|c} a & b & n \\ \alpha & \beta & 0 \end{array} \right) = \sum_{r=0}^{\left\lfloor \frac{n}{2} \right\rfloor -k} \: \lambda_r(0|\alpha, \beta|n) a^r (2a-b)^{\left\lfloor \frac{n}{2} \right\rfloor - r}.
			\end{aligned}
		\end{equation}
	Then, from 	\eqref{ex000}, we get 
		\begin{equation}
		\label{FD2QQ}
		\begin{aligned}
			\Psi(a,b,n)= \sum_{r=0}^{\left\lfloor \frac{n}{2} \right\rfloor -k}  \: \lambda_r(0|\alpha, \beta|n) a^r (2a-b)^{\left\lfloor \frac{n}{2} \right\rfloor - r}.
		\end{aligned}		
\end{equation}		
Consequently, from \eqref{comp3}, we immediately get 		
   \begin{equation}
			\label{Q3}
			\begin{aligned}
				\lambda_r(\:0\:|\alpha, \beta \:|\:n\:) = (-1)^r  \frac{n}{n-r} \binom{n-r}{r}, 
			\end{aligned}
		\end{equation}
which completes the proof.	
			\end{proof}
	
	\begin{theorem}
		\label{H1} 
		For any $n,r, k,\alpha,\beta, n$, the following relation is true

		\begin{equation}
			\label{H2} 
			\begin{aligned}
				\lambda_r(k|\alpha, \beta|n)= (-1)^r \: \frac{\: n \: \: (n-r-k-1)! \: (\left\lfloor \frac{n}{2} \right\rfloor - r )! \: }{\: (n-2r)! \:\: r! \:\:  (\left\lfloor \frac{n}{2} \right\rfloor - r -k)! \: }	\:  \Omega_r\big(k|\:\alpha, \beta \: | n \big).         		
			\end{aligned}
		\end{equation}
	\end{theorem}
	
	\begin{proof}
		To prove \eqref{H2}, define $\bar{\Omega}_r\big(k|\:\alpha, \beta \: | n \big)$ as following
		
		\begin{equation}
			\label{H3} 
			\begin{aligned}
				\lambda_r(k|\alpha, \beta|n)= (-1)^r \: \frac{\: n \: \: (n-r-k-1)! \: (\left\lfloor \frac{n}{2} \right\rfloor - r )! \: }{\: (n-2r)! \:\: r! \:\:  (\left\lfloor \frac{n}{2} \right\rfloor - r -k)! \: }	\:  \bar{\Omega}_r\big(k|\:\alpha, \beta \: | n \big).         		
			\end{aligned}
		\end{equation}
		We need to prove that 
		\begin{equation}
			\label{H4} 
			\begin{aligned}
				\bar{\Omega}_r\big(k|\:\alpha, \beta \: | n \big) \: = \:	\Omega_r\big(k|\:\alpha, \beta \: | n \big),         		
			\end{aligned}
		\end{equation}
		as following. From \eqref{H3}, 	\eqref{FD3}, we get
		
		\begin{equation}
			\label{H5} 
			\begin{aligned}
				& (-1)^r \: \frac{\: n \: \: (n-r-k-1)! \: (\left\lfloor \frac{n}{2} \right\rfloor - r )! \: }{\: (n-2r)! \:\: r! \:\:  (\left\lfloor \frac{n}{2} \right\rfloor - r -k)! \: }	\:  \bar{\Omega}_r\big(k|\:\alpha, \beta \: | n \big) \\   &= \: 
				(2\alpha - \beta) \:
				\Big(\left\lfloor \frac{n}{2} \right\rfloor -k- r +1 \Big)  (-1)^r \: \frac{\: n \: \: (n-r-k)! \: (\left\lfloor \frac{n}{2} \right\rfloor - r )! \: }{\: (n-2r)! \:\: r! \:\:  (\left\lfloor \frac{n}{2} \right\rfloor - r -k +1)! \: }	\:  \bar{\Omega}_r\big(k-1|\:\alpha, \beta \: | n \big) \\
				& + \quad  \alpha \: (r+1) \: 
				(-1)^{r+1} \: \frac{\: n \: \: (n-r-k-1)! \: (\left\lfloor \frac{n}{2} \right\rfloor - r -1 )! \: }{\: (n-2r -2)! \:\: (r+1)! \:\:  (\left\lfloor \frac{n}{2} \right\rfloor - r -k )! \: }	\:  \bar{\Omega}_{r+1}\big(k-1|\:\alpha, \beta \: | n \big).
			\end{aligned}
		\end{equation}	
		Simplifying again, we get 
		\begin{equation}
			\label{H6} 
			\begin{aligned}
				& \:  \bar{\Omega}_r\big(k|\:\alpha, \beta \: | n \big) \\   &= \: 
				(2\alpha - \beta) \: (n-r-k)	\:  \bar{\Omega}_r\big(k-1|\:\alpha, \beta \: | n \big) \\
				& - \quad  \alpha \:
				\: \frac{\: \: (n-2r) \: (n-2r-1) \:}{\: (\left\lfloor \frac{n}{2} \right\rfloor - r ) \: }	\:  \bar{\Omega}_{r+1}\big(k-1|\:\alpha, \beta \: | n \big).
			\end{aligned}
		\end{equation}		
		If $\delta(n)=0$ then $ \delta(n-1)=1$, and if  	$\delta(n)=1$ then $\delta(n-1)=0$. Therefore, for either case, we get the following 
		\begin{equation}
			\label{H7} 
			\begin{aligned}
				\: \frac{\: \: (n-2r) \: (n-2r-1) \:}{\: (\left\lfloor \frac{n}{2} \right\rfloor - r ) \: }	\: &=  \frac{(n-2r-\delta(n)) \:(n-2r-\delta(n-1)) \: }{ \frac{n-\delta(n)}{2}\: - \:r} \\
				&=  2 \:\frac{(n-2r-\delta(n) \:(n-2r-\delta(n-1)) \: }{n-2r-\delta(n) } \\	
				&=   2 \: 	(n-2r-\delta(n-1)).
			\end{aligned}
		\end{equation}
		Hence, from \eqref{H6}, \eqref{H7}, we get
		
		\begin{equation}
			\label{H8} 
			\begin{aligned}
				& \:  \bar{\Omega}_r\big(k|\:\alpha, \beta \: | n \big) \\   &= \: 
				(2\alpha - \beta) \: (n-r-k)	\:  \bar{\Omega}_r\big(k-1|\:\alpha, \beta \: | n \big) \\
				& - \: 2 \:  \alpha \:
				\: 	(n-2r-\delta(n-1))	\:  \bar{\Omega}_{r+1}\big(k-1|\:\alpha, \beta \: | n \big).
			\end{aligned}
		\end{equation}	
		Now, it remains to compute the initial value 
		\[ \bar{\Omega}_r\big(0|\:\alpha, \beta \: | n \big),   \]
		as following.	Put $k=0$ in \eqref{H3}, and noting \eqref{FD3},   we get 
		
		\begin{equation}
			\label{H9} 
			\begin{aligned}
				(-1)^r  \frac{n}{n-r} \binom{n-r}{r}    = (-1)^r \: \frac{\: n \: \: (n-r-1)! \: (\left\lfloor \frac{n}{2} \right\rfloor - r )! \: }{\: (n-2r)! \:\: r! \:\:  (\left\lfloor \frac{n}{2} \right\rfloor - r )! \: }	\:  \bar{\Omega}_r\big( 0 \:|\:\alpha, \beta \: | n \big).         		
			\end{aligned}
		\end{equation}
		Consequently, for any $r$, we get 	
		\begin{equation}
			\label{H10} 
			\begin{aligned}
				1   = \:  \bar{\Omega}_r\big(0\:|\:\alpha, \beta\: | n \big).         		
			\end{aligned}
		\end{equation}
		From \eqref{H10}, \eqref{H8}, and from definition \eqref{omegaDef} of Omega sequence, we immediately conclude
		\begin{equation}
			\label{H11} 
			\begin{aligned}
				\bar{\Omega}_r\big(k|\:\alpha, \beta \: | n \big) \: = \:	\Omega_r\big(k|\:\alpha, \beta\: | n \big).         		
			\end{aligned}
		\end{equation}
		This completes the proof of Theorem \eqref{H1}.	
		\end{proof}
Therefore, from Theorem \eqref{FD1} and Theorem \eqref{H1}, we obtain the following representation for 
\[
\Psi\left( \begin{array}{cc|r} a & b & n \\ \alpha & \beta & k \end{array} \right)
\]
in terms of the \( \Omega \)-sequence, which is quite desirable.

\begin{theorem}[The First Fundamental Theorem of the Quanta Prime Sequence]
	\label{F11} 
	For any numbers \( a, b, \alpha, \beta, n \) such that \( \beta a - \alpha b \neq 0 \), the following expansion holds:
	
	\begin{equation}
		\label{F1100} 
		\begin{aligned}
			\Psi\left( \begin{array}{cc|r} a & b & n \\ \alpha & \beta & k \end{array} \right)  
			&=  \sum_{r=0}^{\lfloor \frac{n}{2} \rfloor - k} (-1)^{r+k} \: \frac{(n-r-k-1)! \: n}{(n-2r)! \: r!} 
			\left(\begin{array}{c} \lfloor \frac{n}{2} \rfloor - r \\ k \end{array}\right) \:
			\Omega_r\big(k \mid \alpha, \beta \mid n\big) \: a^r \: (2a - b)^{\lfloor \frac{n}{2} \rfloor - k - r},
		\end{aligned}
	\end{equation}
	where the coefficients 
	\begin{equation}
		\label{F22} 
		\begin{aligned}
			(-1)^{r+k} \: \frac{(n-r-k-1)! \: n}{(n-2r)! \: r!}
			\left(\begin{array}{c} \lfloor \frac{n}{2} \rfloor - r \\ k \end{array}\right) \:
			\Omega_r\big(k \mid \alpha, \beta \mid n\big)
		\end{aligned}
	\end{equation}
	are integers.
\end{theorem}

\begin{proof}
	Substituting \( \lambda_r(k \mid \alpha, \beta \mid n) \) from equation \eqref{H2} into the formula \eqref{FD2} as
	\[
	(-1)^r \: \frac{n \: (n - r - k - 1)! \: (\lfloor \frac{n}{2} \rfloor - r)!}{(n - 2r)! \: r! \: (\lfloor \frac{n}{2} \rfloor - r - k)!} \: \Omega_r\big(k \mid \alpha, \beta \mid n\big),
	\]
	we immediately obtain the formula \eqref{F1100} in Theorem \ref{F11}, which completes the proof.
\end{proof}

	Now, we obtain the following desirable theorem.
	
	\begin{theorem}[The Second Fundamental Theorem of the Quanta Prime Sequence]
		\label{k0} 
		For any numbers \( \alpha, \beta, n \) with \( (\alpha, \beta) \neq (0, 0) \), the ratio
		\begin{equation}
			\label{k00N} 
			\frac{\Omega_0\left(\left\lfloor \frac{n}{2} \right\rfloor \mid \alpha, \beta \mid n\right)}{(n - 1)(n - 2) \cdots \left(n - \left\lfloor \frac{n}{2} \right\rfloor\right)}
		\end{equation}
		is an integer. Moreover, this ratio gives \( \Psi(\alpha, \beta, n) \). Namely,
		\begin{equation}
			\label{k00} 
			\Psi(\alpha, \beta, n) = \frac{\Omega_0\left(\left\lfloor \frac{n}{2} \right\rfloor \mid \alpha, \beta \mid n\right)}{(n - 1)(n - 2) \cdots \left(n - \left\lfloor \frac{n}{2} \right\rfloor\right)}.
		\end{equation}
	\end{theorem}
	
	\begin{proof}
		To derive the formula \eqref{k00} for \( \Omega_0\left(\left\lfloor \frac{n}{2} \right\rfloor \mid \alpha, \beta \mid n\right) \), we set \( k = \left\lfloor \frac{n}{2} \right\rfloor \) in Theorem \eqref{F11} as follows:
		
		\begin{equation}
			\label{k1} 
			\begin{aligned}
				\Psi\left( \begin{array}{cc|c} a & b & n \\ \alpha & \beta & \left\lfloor \frac{n}{2} \right\rfloor \end{array} \right)  
				&= \sum_{r=0}^{0} (-1)^{r + \left\lfloor \frac{n}{2} \right\rfloor} \: \frac{(n - r - \left\lfloor \frac{n}{2} \right\rfloor - 1)! \: n}{(n - 2r)! \: r!} 
				\left(\begin{array}{c} \lfloor \frac{n}{2} \rfloor - r \\ \left\lfloor \frac{n}{2} \right\rfloor \end{array}\right) 
				\Omega_r\left(\left\lfloor \frac{n}{2} \right\rfloor \mid \alpha, \beta \mid n\right) a^r (2a - b)^{-r} \\
				&= (-1)^{\left\lfloor \frac{n}{2} \right\rfloor} \: \frac{(n - \left\lfloor \frac{n}{2} \right\rfloor - 1)! \: n}{n!} 
				\: \Omega_0\left(\left\lfloor \frac{n}{2} \right\rfloor \mid \alpha, \beta \mid n\right) \\
				&= (-1)^{\left\lfloor \frac{n}{2} \right\rfloor} \: \frac{\Omega_0\left(\left\lfloor \frac{n}{2} \right\rfloor \mid \alpha, \beta \mid n\right)}{(n - 1)(n - 2) \cdots \left(n - \left\lfloor \frac{n}{2} \right\rfloor\right)}.
			\end{aligned}
		\end{equation}
		
		From Theorem \eqref{exp1}, it follows that
		\begin{equation}
			\label{k2} 
			\Psi\left( \begin{array}{cc|c} a & b & n \\ \alpha & \beta & \left\lfloor \frac{n}{2} \right\rfloor \end{array} \right) = (-1)^{\left\lfloor \frac{n}{2} \right\rfloor} \Psi(\alpha, \beta, n).
		\end{equation}
		
		Combining equations \eqref{k1} and \eqref{k2}, we obtain the desired result, completing the proof.
	\end{proof}

\begin{theorem}[Representation for Sums of Powers]
	\label{FA1} 
	\begin{equation}
		\label{FA2} 
		\frac{x^n + y^n}{(x + y)^{\delta(n)}} = \frac{\Omega_0\left(\left\lfloor \frac{n}{2} \right\rfloor \mid xy, -x^2 - y^2 \mid n \right)}{(n - 1)(n - 2) \cdots \left(n - \left\lfloor \frac{n}{2} \right\rfloor \right)}
	\end{equation}
\end{theorem}

\begin{proof}
	From Theorem \eqref{WW3}, we know that
	\[
	\Psi(xy, -x^2 - y^2, n) = \frac{x^n + y^n}{(x + y)^{\delta(n)}}.
	\]
	Setting \( (\alpha, \beta) = (xy, -x^2 - y^2) \) and applying Theorem \eqref{k0}, we immediately obtain the result, completing the proof.
\end{proof}

		Before exploring the second fundamental theorem of Omega sequence, we need to define the Omega space at level $n$ as following.
	\section{The Omega space at level \texorpdfstring{$n$}{}}
	
	\begin{definition}{(The Omega space at Level $n$)} \\
		For any given natural number $n$, we define the Omega space, $\omega(n)$, at level $n$ as following  
		\begin{equation}
			\label{Space} 
			\begin{aligned}
			\omega(n):= \{ (a,b)| \: \Psi(a,b,n) \neq 0 \}.
			\end{aligned}
		\end{equation}
	\end{definition}

Here are some specific examples for points  belonging to Omega space at level $n$. 
\begin{equation}
	\label{S1} 
	\begin{aligned}
		\blacksquare \quad	&\text{If}\: \: &n \equiv \pm 1       \pmod{8} \quad  &\Rightarrow  & 	\quad \: \:\Psi(1,0,n) =1 \neq 0 \:  &\Rightarrow & \quad (1,0) &\in \: \omega(n) \\					
		\blacksquare \quad	&\text{If}\: \:  &n \equiv \pm 2       \pmod{12} \quad  &\Rightarrow  & \quad 	\Psi(1,-1,n) = 1 \neq 0 \:  &\Rightarrow & \quad (1,-1) &\in \: \omega(n)\\
		\blacksquare \quad	&\text{If}\: \:  &n \equiv \pm 3       \pmod{16} \quad  &\Rightarrow  & \quad 	\Psi(1,\sqrt{2},n) = -1-\sqrt{2}\neq 0 \:  &\Rightarrow & \quad (1,\sqrt{2}) &\in \: \omega(n)\\
		\blacksquare \quad	&\text{If}\: \:  &n \equiv \pm 4      \pmod{20} \quad  &\Rightarrow  & \quad 	\Psi(1,\varphi -1,n) = -\varphi \neq 0 \:  &\Rightarrow & \quad (1,\varphi -1) &\in \: \omega(n)\\	
		\blacksquare \quad	&\text{If}\: \:  &n \equiv \pm 5       \pmod{24} \quad  &\Rightarrow  & \quad 	\Psi(1,\sqrt{3},n) =  2+ \sqrt{3}\neq 0 \:  &\Rightarrow & \quad (1,\sqrt{3}) &\in \: \omega(n)
	\end{aligned}
\end{equation}
where $\varphi$ is the Golden ratio.

Here are some specific examples for points not belonging to Omega space at level $n$. 
\begin{equation}
	\label{S11} 
	\begin{aligned}
		\blacksquare \quad	&\text{If}\: \: &n \equiv \pm 2       \pmod{8} \quad  &\Rightarrow  & 	\quad \: \:\Psi(1,0,n) =0 \:  &\Rightarrow & \quad (1,0) &\notin \: \omega(n) \\					
		\blacksquare \quad	&\text{If}\: \:  &n \equiv \pm 3       \pmod{12} \quad  &\Rightarrow  & \quad 	\Psi(1,-1,n) =0 \:  &\Rightarrow & \quad (1,-1) &\notin \: \omega(n)\\
		\blacksquare \quad	&\text{If}\: \:  &n \equiv \pm 4       \pmod{16} \quad  &\Rightarrow  & \quad 	\Psi(1,\sqrt{2},n) =0 \:  &\Rightarrow & \quad (1,\sqrt{2}) &\notin \: \omega(n)\\
		\blacksquare \quad	&\text{If}\: \:  &n \equiv \pm 5       \pmod{20} \quad  &\Rightarrow  & \quad 	\Psi(1,\varphi -1,n) =0 \:  &\Rightarrow & \quad (1,\varphi -1) &\notin \: \omega(n)\\	
		\blacksquare \quad	&\text{If}\: \:  &n \equiv \pm 6       \pmod{24} \quad  &\Rightarrow  & \quad 	\Psi(1,\sqrt{3},n) =0 \:  &\Rightarrow & \quad (1,\sqrt{3}) &\notin \: \omega(n)
	\end{aligned}
\end{equation}
 
These examples suggest one to define the Kernel of Omega space for further research developments. 
\begin{definition}{(The Kernel of the Omega space at level $n$)} 
	For any given natural number $n$, we define the Kernel of the Omega space, $Ker_{\omega}(n)$, at level $n$ as following  
	\begin{equation}
		\label{kernel} 
		\begin{aligned}
			Ker_{\omega}(n):= \{ (a,b)| \: \Psi(a,b,n) = 0 \}.
		\end{aligned}
	\end{equation}
\end{definition}

\begin{theorem}
	\label{sH} 
	For any natural number \( n \), the space \( \omega(n) \) is infinite.
\end{theorem} 

\begin{proof}
	By \eqref{WW4}, we know that
	\[
	\Psi(xy, -x^2 - y^2, n) \neq 0
	\]
	for any integers \( x \) and \( y \) such that \( x \neq -y \). This implies that the space \( \omega(n) \) contains all integer points \( (xy, -x^2 - y^2) \) satisfying \( x \neq -y \), thus proving that \( \omega(n) \) is infinite.
\end{proof}

\begin{theorem}[The Second Fundamental Theorem of the Quanta Prime Sequence (Version 2)]
	\label{space3} 
	For any natural number \( n \) and any point \( (\alpha, \beta) \in \omega(2n) \), we have
	\begin{equation}
		\label{space4} 
		n(n+1) \dotsm (2n-1) = \frac{\Omega_0\left(n \mid \alpha, \beta \mid 2n\right)}{\Psi(\alpha, \beta, 2n)}.
	\end{equation}    
\end{theorem}

\begin{proof}
	For any natural number \( n > 1 \), we observe from formula \eqref{k00} that
	\[
	\Psi(\alpha, \beta, n) \neq 0 \quad \iff \quad \Omega_0\left(\left\lfloor \frac{n}{2} \right\rfloor \mid \alpha, \beta \mid n\right) \neq 0.
	\]
	Thus, the following ratio 
	\begin{equation}
		\label{space1} 
		\frac{\Omega_0\left(\left\lfloor \frac{n}{2} \right\rfloor \mid \alpha, \beta \mid n\right)}{\Psi(\alpha, \beta, n)}
	\end{equation}    
	is well-defined for any point \( (\alpha, \beta) \in \omega(n) \). Consequently, for any \( (\alpha, \beta) \in \omega(n) \), the product 
	\[
	(n - 1)(n - 2) \dotsm \left(n - \left\lfloor \frac{n}{2} \right\rfloor\right)
	\]
	can be represented as follows:
	
	\begin{equation}
		\label{space2} 
		(n - 1)(n - 2) \dotsm \left(n - \left\lfloor \frac{n}{2} \right\rfloor\right) = \frac{\Omega_0\left(\left\lfloor \frac{n}{2} \right\rfloor \mid \alpha, \beta \mid n\right)}{\Psi(\alpha, \beta, n)}.
	\end{equation}    
	
	Substituting \( n \) with \( 2n \) in \eqref{space2} completes the proof of the theorem.
\end{proof}

The following unexpected theorem proves that any new prime must be a factor to Quanta Prime sequence.
		
	\section{The k-th Prime Number}

	The \( k \)-th prime number, denoted as \( p_k \), refers to the \( k \)-th element in the ordered sequence of prime numbers. A prime number is defined as a natural number greater than 1 that has no positive divisors other than 1 and itself \cite{hardy2008mathematical}. The sequence of prime numbers begins as follows:
	\[
	2, 3, 5, 7, 11, 13, 17, 19, 23, 29, \ldots
	\]
For example, the first few values of \( p_k \) are: \( p_1 = 2 \), \( p_2 = 3 \), \( p_3 = 5 \), and \( p_4 = 7 \).
	
	The determination of the \( k \)-th prime can be approached using various algorithms, including the Sieve of Eratosthenes, which efficiently generates all prime numbers up to a specified limit \cite{sieve2019}. Additionally, the distribution of prime numbers has been extensively studied, with the Prime Number Theorem stating that the number of primes less than a given integer \( n \) is asymptotically approximated by \( \frac{n}{\log n} \) \cite{rosser1938estimates}. This foundational result provides significant insights into the behavior of prime numbers as \( k \) increases.

\section{The Emergence of a New Prime Number}
Prime numbers have consistently fascinated mathematicians, not only due to their fundamental importance in number theory but also because of their enigmatic and seemingly irregular distribution among the integers. In this section, we present a noteworthy result within the Quanta Prime Sequence, which illuminates a structured mechanism for the emergence of primes. Specifically, this result reveals a distinctive property:

\begin{center}
	For each \( k \)-th prime \( p_k \), the subsequent prime \( p_{k+1} \) divides the Quanta Prime Sequence associated with \( p_k \). Additionally, an infinite set of parameter pairs \( (\alpha, \beta) \) within the sequence allows for \textit{considerable flexibility in the analysis of prime emergence}. The general proof of this property is presented in Theorem \eqref{infinite_params}, with a specialized version given in Theorem \eqref{gen2} on the following page.
\end{center}

This divisibility condition underscores the intrinsic link between the Quanta Prime Sequence and the progression of primes, offering a perspective on how primes may systematically 'emerge' within a defined numerical framework. The infinite parameter space not only highlights the depth of this connection but also provides a powerful means to observe prime patterns and periodicities with fine-grained control. The presence of an infinite number of configurations for \( (\alpha, \beta) \) opens up extensive opportunities for researchers to examine recursive relationships that give rise to prime numbers and to test hypotheses across a wide range of conditions within the sequence. By framing the appearance of new primes within the parameters of the Quanta Prime Sequence, we gain a valuable tool for probing the underlying dynamics of prime emergence. This structured framework invites researchers to delve into the mechanics of prime formation and to uncover deeper insights into the distribution of these essential building blocks of number theory. Prime numbers, due to their fundamental importance, have led to extensive studies of their distribution in number theory. One foundational result in this area is the Bertrand–Chebyshev Theorem, which asserts the existence of at least one prime in the interval \( (n, 2n) \) for any integer \( n > 1 \). Originally conjectured by Joseph Bertrand in 1845 and later rigorously proved by Pafnuty Chebyshev in 1852, this theorem (see \cite{24}, \cite{chebyshev1852}, \cite{apostol1976}, \cite{nathanson2000}) provides a key insight into prime intervals and supports the investigation of prime behavior within sequences like the Quanta Prime Sequence.

\begin{theorem}[Bertrand–Chebyshev Theorem]
	For any integer \( n > 1 \), there exists at least one prime \( p \) such that
	\[
	n < p < 2n.
	\]
\end{theorem}

\begin{theorem}[Infinite Parameter Space for Prime Emergence]
	\label{infinite_params}
	For any \( k \)-th prime \( p_k \), there exists an infinite space \( \varpi \) of points \( (\alpha, \beta) \) and an index set \( I \) such that the following property holds:
	\begin{equation} 
		\label{space8} 
		p_{k+1} \;\; \bigg| \, \sum_{\substack{(\alpha, \beta) \in \varpi \\ \lambda \in I}} \lambda \, \frac{\Omega_0\big(p_k \mid \alpha, \beta \mid 2p_k \big)}{\Psi(\alpha, \beta, 2p_k)}.
	\end{equation}    
\end{theorem}

This theorem vastly extends classical results by introducing an infinite set of parameters \( (\alpha, \beta) \in \varpi \) that satisfy the divisibility property, thus offering an unlimited range of configurations for exploring prime emergence within the Quanta Prime Sequence. The existence of an infinite parameter space allows researchers to not only investigate the emergence of a prime \( p_{k+1} \) with unprecedented control but also to examine complex patterns and behaviors of primes across diverse intervals and conditions. This flexibility opens new avenues for studying prime generation mechanisms, offering deep insights into number theory and prime distribution beyond the constraints of classical theorems. This result ensures a minimum density of primes as numbers grow larger, and it has inspired further research into prime gaps and distribution (\cite{nathanson2000}, \cite{toth2019}).

\begin{proof}
	From the Bertrand–Chebyshev theorem, we have
	\begin{equation}
		\label{space44e} 
		p_{k+1} \;\; \, \big| \;\;\, p_k \; (p_k + 1)\; (p_k + 2) \dotsm (2p_k - 1).
	\end{equation}
	
For any point \( (\alpha, \beta) \in \omega(2p_k) \), the ratio
\begin{equation}
	\label{space6E} 
	\frac{\Omega_0\big(p_k \mid \alpha, \beta \mid 2p_k \big)}{\Psi(\alpha, \beta, 2p_k)}
\end{equation}
is an integer. 

Setting \( n = p_k \) in \eqref{space4}, we obtain
\begin{equation}
	\label{space444} 
	p_k(p_k + 1) \dotsm (2p_k - 1) = \frac{\Omega_0\big(p_k \mid \alpha, \beta \mid 2p_k \big)}{\Psi(\alpha, \beta, 2p_k)},
\end{equation}    
where \( (\alpha, \beta) \in \omega(p_k) \).

Thus, combining \eqref{space44e} and \eqref{space444}, we get
\begin{equation}
	\label{space6} 
	p_{k+1} \;\;\; \, \bigg| \, \; \frac{\Omega_0\big(p_k \mid \alpha, \beta \mid 2p_k \big)}{\Psi(\alpha, \beta, 2p_k)}.
\end{equation}    
Moreover, for any finite set \( \varpi \subset \omega(2p_k) \) and any finite set \( I \) of integers, we have
\begin{equation}
	\label{space8} 
	p_{k+1} \;\;\, \bigg| \, \sum_{\substack{(\alpha, \beta) \in \varpi \\ \lambda \in I}} \lambda \, \frac{\Omega_0\big(p_k \mid \alpha, \beta \mid 2p_k \big)}{\Psi(\alpha, \beta, 2p_k)}
\end{equation}    	
	  which completes the proof.
	\end{proof}

We should also observe the following Theorem 

\begin{theorem}
	\label{gen1} 
	For any \( k > 1 \) and any point \( (\alpha, \beta) \in \omega(2p_k) \), the ratio 
	\begin{equation}
		\label{gen11} 
		\frac{\Omega_0\big(p_k \mid \alpha, \beta \mid 2p_k \big)}{p_k (2p_k - 1)(2p_k - 2) \, \Psi(\alpha, \beta, 2p_k)}
	\end{equation}
	is an integer. Moreover,
	\begin{equation}
		\label{gen111} 
		p_{k+1} \;\;\, \big|\;\; \, \frac{\Omega_0\big(p_k \mid \alpha, \beta \mid 2p_k \big)}{p_k (2p_k - 1)(2p_k - 2) \, \Psi(\alpha, \beta, 2p_k)}.
	\end{equation}    
\end{theorem} 

\begin{proof}
	Since \( (p_k, (2p_k - 1)(2p_k - 2)) = 1 \), the result follows immediately.
\end{proof}

\begin{theorem}[Enhanced Prime Emergence with Two Parameters]
	\label{gen2}
	For any point \( (\alpha, \beta) \), the following property holds:
	\begin{equation}
		\label{gen22} 
		p_{k+1}\;\; \, \big| \;\;\, \Omega_0\big(p_k \mid \alpha, \beta \mid 2p_k \big).
	\end{equation}
	This theorem extends the Bertrand–Chebyshev result by introducing two parameters, \( (\alpha, \beta) \), that provide a unique degree of flexibility in examining the conditions under which new primes emerge. By adjusting these parameters, researchers can pinpoint instances of prime emergence within the Quanta Prime Sequence, offering a powerful framework for exploring prime behavior across specified intervals. This structured approach allows for a more nuanced analysis of prime distribution and periodicities, adding a layer of control and insight that surpasses classical results.
\end{theorem}

\begin{proof}
	From Theorem \eqref{infinite_params}, we know that
	\begin{equation}
		p_{k+1} \, \big| \, \frac{\Omega_0\big(p_k \mid \alpha, \beta \mid 2p_k \big)}{\Psi(\alpha, \beta, 2p_k)}.
	\end{equation}
	Since it also holds that
	\begin{equation}
		\frac{\Omega_0\big(p_k \mid \alpha, \beta \mid 2p_k \big)}{\Psi(\alpha, \beta, 2p_k)} \, \big| \, \Omega_0\big(p_k \mid \alpha, \beta \mid 2p_k \big),
	\end{equation}
	it follows that \( p_{k+1} \, \big| \, \Omega_0\big(p_k \mid \alpha, \beta \mid 2p_k \big) \), thereby completing the proof.
\end{proof}

\section{\textbf{Fascinating Patterns in the Quanta Prime Sequence}}

In this section, we explore a collection of intriguing and distinct cases that emerge from the Second Fundamental Theorem of the Quanta Prime Sequence (Theorem \ref{k0}). Each special case provides a unique perspective, uncovering patterns and properties that deepen our understanding of the sequence’s underlying structure. These examples demonstrate the theorem’s adaptability and its capacity to reveal unexpected relationships. Through these cases, we aim to shed light on the broader implications of the Quanta Prime Sequence, offering insights that may inspire further research and innovative applications.

\subsection{\textbf{The Quanta Prime Sequence Associated with the Point (1, 1)}}

\begin{theorem}
	\label{PP00}
	For any natural number \( n \), the Quanta Prime Sequence satisfies:
	\begin{equation}
		\label{PP1} 
		\frac{\Omega_0\left(\left\lfloor \frac{n}{2} \right\rfloor \mid 1, 1 \mid n \right)}{(n-1)(n-2) \dotsm \left(n - \left\lfloor \frac{n}{2} \right\rfloor\right)} = 
		\begin{cases}
			+2, & n \equiv 0 \pmod{6}, \\
			+1, & n \equiv \pm 1 \pmod{6}, \\
			-1, & n \equiv \pm 2 \pmod{6}, \\
			-2, & n \equiv \pm 3 \pmod{6}.
		\end{cases}
	\end{equation}
\end{theorem}

\begin{proof}
	We utilize the following established property of the \( \Psi \)-sequence:
	\begin{equation}
		\label{PP0}
		\Psi(1,1,n) =
		\begin{cases}
			+2, & n \equiv 0 \pmod{6}, \\
			+1, & n \equiv \pm 1 \pmod{6}, \\
			-1, & n \equiv \pm 2 \pmod{6}, \\
			-2, & n \equiv \pm 3 \pmod{6}.
		\end{cases}   
	\end{equation}
	Substituting \( \alpha = 1 \) and \( \beta = 1 \) into formula \eqref{k00} in Theorem \eqref{k0} confirms the stated result.
\end{proof}

\subsection{\textbf{The Quanta Prime Sequence Associated with the Point (1, 0)}} \leavevmode\\

\vspace{10pt}

In a similar fashion, we obtain the following result.
\begin{theorem}
	\label{PP00Q}
	For any natural number \( n \), the Quanta Prime Sequence satisfies:
	\begin{equation}
		\label{PP1Q} 
		\frac{\Omega_0\left(\left\lfloor \frac{n}{2} \right\rfloor \mid 1, 0 \mid n \right)}{(n-1)(n-2) \dotsm \left(n - \left\lfloor \frac{n}{2} \right\rfloor\right)} = 
		\begin{cases}
			+2, & n \equiv 0 \pmod{8}, \\
			+1, & n \equiv \pm 1 \pmod{8}, \\
			\;0, & n \equiv \pm 2 \pmod{8}, \\				
			-1, & n \equiv \pm 3 \pmod{8}, \\
			-2, & n \equiv \pm 4 \pmod{8}.
		\end{cases}        
	\end{equation}
\end{theorem}

\begin{proof}
	The following known property of the \( \Psi \)-sequence is relevant:
	\begin{equation}
		\label{PP0Q}
		\Psi(1,0,n) =
		\begin{cases}
			+2, & n \equiv 0 \pmod{8}, \\
			+1, & n \equiv \pm 1 \pmod{8}, \\
			\;0, & n \equiv \pm 2 \pmod{8}, \\				
			-1, & n \equiv \pm 3 \pmod{8}, \\
			-2, & n \equiv \pm 4 \pmod{8}.
		\end{cases}     
	\end{equation}
	Setting \( \alpha = 1 \) and \( \beta = 0 \) in formula \eqref{k00} in Theorem \eqref{k0} confirms the result.
\end{proof}

\subsection{\textbf{The Quanta Prime Sequence Associated with the Point (1, -1)}}
\leavevmode\\
		\begin{theorem}
		\label{PP}
	For any natural number \( n \), the Quanta Prime Sequence satisfies:
		\begin{equation}
			\label{PP1A} 
			\frac{\Omega_0\left(\left\lfloor \frac{n}{2} \right\rfloor \mid 1, -1 \mid n \right)}{(n-1)(n-2) \dotsm \left(n - \left\lfloor \frac{n}{2} \right\rfloor\right)} = 
			\begin{cases}
				+2, & n \equiv 0 \pmod{12}, \\
				+1, & n \equiv \pm 1, \pm 2 \pmod{12}, \\
				\;\: 0, & n \equiv \pm 3 \pmod{12}, \\
				-1, & n \equiv \pm 4, \pm 5 \pmod{12}, \\
				-2, & n \equiv \pm 6 \pmod{12}.
			\end{cases}  
		\end{equation}
	\end{theorem}
	
	\begin{proof}
		This result follows from the periodic behavior of the \( \Psi \)-sequence, specifically:
		\begin{equation}
			\label{peroidicity1}
			\Psi(1, -1, n) = 
			\begin{cases}
				+2, & n \equiv 0 \pmod{12}, \\
				+1, & n \equiv \pm 1, \pm 2 \pmod{12}, \\
				\;\: 0, & n \equiv \pm 3 \pmod{12}, \\
				-1, & n \equiv \pm 4, \pm 5 \pmod{12}, \\
				-2, & n \equiv \pm 6 \pmod{12}.
			\end{cases}
		\end{equation}
		Substituting \( \alpha = 1 \) and \( \beta = -1 \) into formula \eqref{k00} from Theorem \eqref{k0} verifies the result.
	\end{proof}

	\subsection{\textbf{The Quanta Prime Sequence Associated with the Point $(1, \sqrt{2})$}}
	\leavevmode\\
			\begin{theorem}
		\label{PP}
	For any natural number \( n \), the Quanta Prime Sequence satisfies:
		\begin{equation}
			\label{PP1A} 
			\frac{\Omega_0\left(\left\lfloor \frac{n}{2} \right\rfloor \mid 1, \sqrt{2} \mid n \right)}{(n-1)(n-2) \dotsm \left(n - \left\lfloor \frac{n}{2} \right\rfloor\right)} = 
			\begin{cases}
				2 & n \equiv 0 \pmod{16}, \\
				1 & n \equiv \pm 1 \pmod{16}, \\
				-\sqrt{2} & n \equiv \pm 2 \pmod{16}, \\
				-1 - \sqrt{2} & n \equiv \pm 3 \pmod{16}, \\
				0 & n \equiv \pm 4 \pmod{16}, \\
				1 + \sqrt{2} & n \equiv \pm 5 \pmod{16}, \\
				\sqrt{2} & n \equiv \pm 6 \pmod{16}, \\
				-1 & n \equiv \pm 7 \pmod{16}, \\
				-2 & n \equiv \pm 8 \pmod{16}.
			\end{cases}
		\end{equation}
	\end{theorem}
	
	\begin{proof}
		This result relies on the established periodic properties of the \( \Psi \)-sequence for the values \( (1, \sqrt{2}) \), given by
		\begin{equation}
			\label{root-2}
			\Psi(1, \sqrt{2}, n) = 
			\begin{cases}
				2 & n \equiv 0 \pmod{16}, \\
				1 & n \equiv \pm 1 \pmod{16}, \\
				-\sqrt{2} & n \equiv \pm 2 \pmod{16}, \\
				-1 - \sqrt{2} & n \equiv \pm 3 \pmod{16}, \\
				0 & n \equiv \pm 4 \pmod{16}, \\
				1 + \sqrt{2} & n \equiv \pm 5 \pmod{16}, \\
				\sqrt{2} & n \equiv \pm 6 \pmod{16}, \\
				-1 & n \equiv \pm 7 \pmod{16}, \\
				-2 & n \equiv \pm 8 \pmod{16}.
			\end{cases}
		\end{equation}
		Substituting \( \alpha = 1 \) and \( \beta = \sqrt{2} \) into formula \eqref{k00} in Theorem \eqref{k0} leads directly to the stated result.
	\end{proof}

\subsection{\textbf{Revealing Golden Ratio Patterns in $\Psi(1, \phi - 1, n)$}} \leavevmode\\
The golden ratio, \( \phi = \frac{1 + \sqrt{5}}{2} \), has captivated mathematicians, artists, and scientists for centuries with its unique properties and recurring appearance in nature, art, and architecture. Here, we uncover how this fascinating number shapes the behavior of the Quanta Prime Sequence, creating a distinctive alternating pattern in the expression \( \Psi(1, \phi - 1, n) \).

\begin{theorem}
	\label{PP}
	For any natural number \( n \), the Quanta Prime Sequence satisfies:
	\begin{equation}
		\label{PP1AA} 
		\frac{\Omega_0\left(\left\lfloor \frac{n}{2} \right\rfloor \mid 1, \phi - 1 \mid n \right)}{(n-1)(n-2) \dotsm \left(n - \left\lfloor \frac{n}{2} \right\rfloor\right)} = 
		\begin{cases}
			2 & n \equiv 0 \pmod{20}, \\
			1 & n \equiv \pm 1 \pmod{20}, \\
			-\phi + 1 & n \equiv \pm 2 \pmod{20}, \\
			-\phi & n \equiv \pm 3, \pm 4 \pmod{20}, \\
			0 & n \equiv \pm 5 \pmod{20}, \\
			\phi & n \equiv \pm 6, \pm 7 \pmod{20}, \\
			\phi - 1 & n \equiv \pm 8 \pmod{20}, \\
			-1 & n \equiv \pm 9 \pmod{20}, \\
			-2 & n \equiv \pm 10 \pmod{20}.
		\end{cases}
	\end{equation}
\end{theorem}

\begin{proof}
	This result is derived from the periodic properties inherent to the \( \Psi \)-sequence when evaluated with parameters \( (1, \phi - 1) \). The periodic behavior is observed as follows:
	\begin{equation}
		\label{phi}
		\Psi(1, \phi - 1, n) =
		\begin{cases}
			2 & n \equiv 0 \pmod{20}, \\
			1 & n \equiv \pm 1 \pmod{20}, \\
			-\phi + 1 & n \equiv \pm 2 \pmod{20}, \\
			-\phi & n \equiv \pm 3, \pm 4 \pmod{20}, \\
			0 & n \equiv \pm 5 \pmod{20}, \\
			\phi & n \equiv \pm 6, \pm 7 \pmod{20}, \\
			\phi - 1 & n \equiv \pm 8 \pmod{20}, \\
			-1 & n \equiv \pm 9 \pmod{20}, \\
			-2 & n \equiv \pm 10 \pmod{20}.
		\end{cases}
	\end{equation}
	Substituting \( \alpha = 1 \) and \( \beta = \phi - 1 \) into formula \eqref{k00} from Theorem \eqref{k0} directly leads to the stated result.
\end{proof}
	\subsection{\textbf{The Quanta Prime Sequence Associated with the Point $(1, \sqrt{3})$}}
	\leavevmode\\
	The function \( \Psi(1, \sqrt{3}, n) \) exhibits a notable periodic behavior with a cycle of 24. This periodicity introduces distinct alternating patterns, each dependent on the values of \( n \), as demonstrated in the following theorem.

\begin{theorem}
	\label{PP}
	For any natural number \( n \), the Quanta Prime Sequence satisfies:
	\begin{equation}
		\label{PP1AAA} 
		\frac{\Omega_0\left(\left\lfloor \frac{n}{2} \right\rfloor \mid 1, \sqrt{3} \mid n \right)}{(n-1)(n-2) \dotsm \left(n - \left\lfloor \frac{n}{2} \right\rfloor\right)} = 
		\begin{cases}
			2 & n \equiv 0 \pmod{24}, \\
			1 & n \equiv \pm 1, \pm 4 \pmod{24}, \\
			-\sqrt{3} & n \equiv \pm 2 \pmod{24}, \\
			-1 - \sqrt{3} & n \equiv \pm 3 \pmod{24}, \\
			2 + \sqrt{3} & n \equiv \pm 5 \pmod{24}, \\
			0 & n \equiv \pm 6 \pmod{24}, \\
			-2 - \sqrt{3} & n \equiv \pm 7 \pmod{24}, \\
			-1 & n \equiv \pm 8, \pm 11 \pmod{24}, \\
			1 + \sqrt{3} & n \equiv \pm 9 \pmod{24}, \\
			\sqrt{3} & n \equiv \pm 10 \pmod{24}, \\
			-2 & n \equiv \pm 12 \pmod{24}.
		\end{cases}
	\end{equation}
\end{theorem}

\begin{proof}
	This result follows from the periodic nature of the Quanta Prime Sequence, specifically when evaluated with the parameters \( (1, \sqrt{3}) \). The periodic behavior can be expressed as:
	\begin{equation}
		\label{root-3}
		\Psi(1, \sqrt{3}, n) = 
		\begin{cases}
			2 & n \equiv 0 \pmod{24}, \\
			1 & n \equiv \pm 1, \pm 4 \pmod{24}, \\
			-\sqrt{3} & n \equiv \pm 2 \pmod{24}, \\
			-1 - \sqrt{3} & n \equiv \pm 3 \pmod{24}, \\
			2 + \sqrt{3} & n \equiv \pm 5 \pmod{24}, \\
			0 & n \equiv \pm 6 \pmod{24}, \\
			-2 - \sqrt{3} & n \equiv \pm 7 \pmod{24}, \\
			-1 & n \equiv \pm 8, \pm 11 \pmod{24}, \\
			1 + \sqrt{3} & n \equiv \pm 9 \pmod{24}, \\
			\sqrt{3} & n \equiv \pm 10 \pmod{24}, \\
			-2 & n \equiv \pm 12 \pmod{24}.
		\end{cases}
	\end{equation}
	Substituting \( \alpha = 1 \) and \( \beta = \sqrt{3} \) into formula \eqref{k00} from Theorem \eqref{k0} confirms the periodic result as stated.
\end{proof}

\subsection{\textbf{Uncovering Periodic Patterns with the Fibonacci and Lucas Sequences}}

\begin{theorem}{(Representation of a Combination of Fibonacci and Lucas Sequences)}
	\label{FL} 
	For any natural number \( n \), the Quanta Prime Sequence satisfies:
	\begin{equation}
		\label{FL1} 
		\frac{\Omega_0\left(\left\lfloor \frac{n}{2} \right\rfloor \mid 1, \sqrt{5} \mid n\right)}{(n-1)(n-2) \dotsm \left(n - \left\lfloor \frac{n}{2} \right\rfloor\right)} = 
		\begin{cases}
			L\left(\frac{n}{2}\right) & n \equiv 0 \pmod{4}, \\
			L\left(\frac{n+1}{2}\right) + F\left(\frac{n-1}{2}\right) \sqrt{5} & n \equiv 1 \pmod{4}, \\
			-F\left(\frac{n}{2}\right) \sqrt{5} & n \equiv 2 \pmod{4}, \\
			-L\left(\frac{n-1}{2}\right) - F\left(\frac{n+1}{2}\right) \sqrt{5} & n \equiv 3 \pmod{4}.
		\end{cases}
	\end{equation}
\end{theorem}

\begin{proof}
	This result derives from the periodic properties of the Quanta Prime Sequence, particularly when evaluated with parameters \( (1, \sqrt{5}) \). The periodic behavior is given by:
	\begin{equation}
		\label{(1,Root5)}
		\Psi(1, \sqrt{5}, n) = 
		\begin{cases}
			L\left(\frac{n}{2}\right) & n \equiv 0 \pmod{4}, \\
			L\left(\frac{n+1}{2}\right) + F\left(\frac{n-1}{2}\right) \sqrt{5} & n \equiv 1 \pmod{4}, \\
			-F\left(\frac{n}{2}\right) \sqrt{5} & n \equiv 2 \pmod{4}, \\
			-L\left(\frac{n-1}{2}\right) - F\left(\frac{n+1}{2}\right) \sqrt{5} & n \equiv 3 \pmod{4}.
		\end{cases}
	\end{equation}
	By substituting \( \alpha = 1 \) and \( \beta = \sqrt{5} \) into formula \eqref{k00} from Theorem \eqref{k0}, we confirm the periodic result as stated.
\end{proof}

\begin{figure}[h!]
	\centering
	\includegraphics[width=0.8\textwidth]{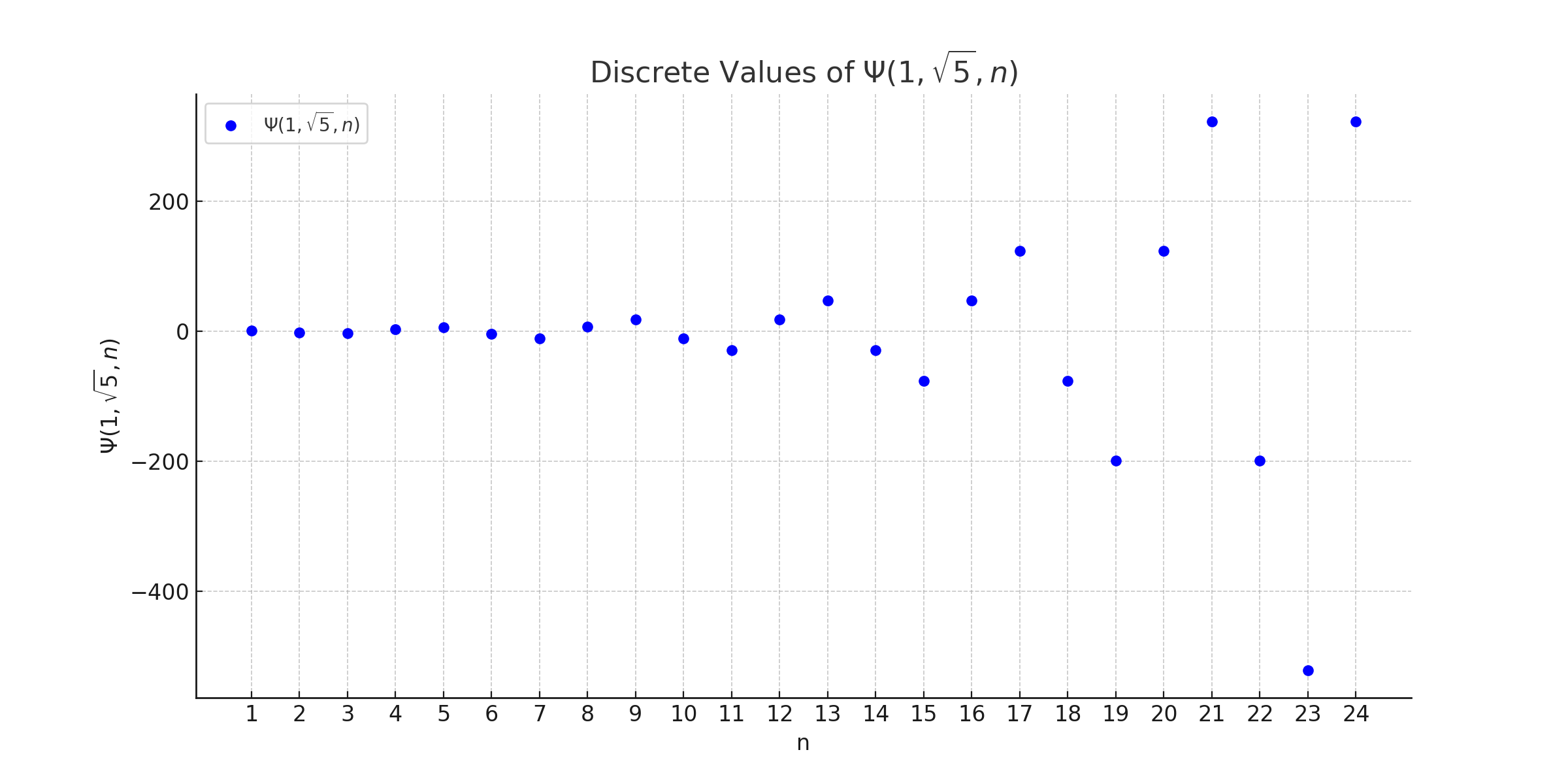}
	\caption{Discrete values of $\Psi(1, \sqrt{5}, n)$ for natural values of $n$.}
	\label{fig:Psi_discrete}
\end{figure}

\subsection{\textbf{The Quanta Prime Sequence Associated with the Point \texorpdfstring{$(1,2)$}{}}}
\leavevmode\\

\begin{theorem}
	\label{DA}
	For any natural number \( n \), the Quanta Prime Sequence satisfies:
	\begin{equation}
		\label{AU10} 
		\frac{\Omega_0\left(\left\lfloor \frac{n}{2} \right\rfloor \mid 1, 2 \mid n \right)}{(n-1)(n-2) \dotsm \left(n - \left\lfloor \frac{n}{2} \right\rfloor\right)} = (-1)^{\left\lfloor \frac{n}{2} \right\rfloor} \cdot 2^{\delta(n-1)} \cdot n^{\delta(n)}.
	\end{equation}
\end{theorem}	

\begin{proof}
	From equation \eqref{comp3}, we have:
	\[
	\Psi(1,2,n) = \frac{n}{n - \left\lfloor \frac{n}{2} \right\rfloor} \binom{n - \left\lfloor \frac{n}{2} \right\rfloor}{\left\lfloor \frac{n}{2} \right\rfloor} \cdot (-1)^{\left\lfloor \frac{n}{2} \right\rfloor}.
	\]
	
	After simplification, this reduces to:
	\[
	\Psi(1,2,n) = (-1)^{\left\lfloor \frac{n}{2} \right\rfloor} \cdot 2^{\delta(n-1)} \cdot n^{\delta(n)},
	\]
	where
	\[
	\delta(n) = \begin{cases} 
		0 & \text{for \( n \) even}, \\ 
		1 & \text{for \( n \) odd}.
	\end{cases}
	\]
	
	Substituting \( \alpha = 1 \) and \( \beta = 2 \) into formula \eqref{k00} from Theorem \eqref{k0} confirms the result as stated.
\end{proof}

Theorem \eqref{DA} is, in fact, a specific instance of the more general result presented in the following theorem.

\begin{theorem}
	\label{AU9}
	The Omega sequence associated with \( n \) and the point \( (1,2) \) is given by:
	\begin{equation}
		\label{AU10A} 
		\Omega_r\left(k \mid 1, 2 \mid n\right) = (-2)^k \prod_{\lambda=0}^{k-1} \left(n - \delta(n+1) - 2r - 2 \lambda \right).
	\end{equation}
\end{theorem}

\begin{proof}
	For the point \( (1,2) \), setting \( 2\eta - \zeta = 0 \) in the definition \eqref{G0}, we observe that the Omega Prime sequence satisfies the recurrence:
	\begin{equation}
		\label{AU1W} 
		\begin{aligned}
			\Omega_r\left(k \mid 1, 2 \mid n\right) &= -2 \left(n - 2r - \delta(n-1)\right) \Omega_{r+1}\left(k-1 \mid 1, 2 \mid n\right), \\
			\Omega_r\left(0 \mid 1, 2 \mid n\right) &= 1 \quad \text{for all} \quad r.
		\end{aligned}
	\end{equation}
	
	Solving \eqref{AU1W} iteratively, we obtain the formula stated in \eqref{AU10A}, as required.
\end{proof}

\subsection{The Quanta Prime Sequence Associated with the Point \texorpdfstring{$(0, -1)$}{}}
\leavevmode\\
\begin{theorem}
	\label{AU11}
	The Omega sequence associated with \( n \) and \( (0, -1) \) is given by:
	\begin{equation}
		\label{AU12} 
		\Omega_r\big(k \mid 0, -1 \mid n\big) = \prod_{\lambda = 1}^{k} \left(n - r - \lambda \right).
	\end{equation}
	
	Additionally,
	\begin{equation}
		\label{AU13} 
		\Omega_0\left(\left\lfloor \frac{n}{2} \right\rfloor \mid 0, -1 \mid n\right) = (n-1)(n-2) \cdots \left(n - \left\lfloor \frac{n}{2} \right\rfloor\right).
	\end{equation}
\end{theorem}
\begin{proof}
	By substituting \( \alpha = 0 \) and \( \beta = -1 \) into formula \eqref{k00} from Theorem \eqref{k0}, and observing that \( \Psi(0, -1, n) = 1 \) for any natural number \( n \), we immediately obtain the stated result.
\end{proof}

	\section{A New Classical Combinatorial Identity}	
	
	\begin{theorem}[New Combinatorial Identity]
		\label{AU7}	
		\begin{equation}
			\label{AU8} 	
			\begin{aligned}	
				(n-1)(n-2) \cdots \left(n - \left\lfloor \frac{n}{2} \right\rfloor\right) &= \\ 
				2^{\left\lfloor \frac{n}{2} \right\rfloor - \delta(n+1)} \left(n + \delta(n-1) - 2\right) 
				&\left(n + \delta(n-1) - 4\right) \left(n + \delta(n-1) - 6\right) \cdots (5)(3)(1).
			\end{aligned}	
		\end{equation}
	\end{theorem}
	
	\begin{proof}
		By substituting from \eqref{def0}, we have:
		\begin{equation} \label{ABAB}
			\Psi(1, -2, n) = 2^{\delta(n+1)}.
		\end{equation}
		
		Utilizing \eqref{ABAB} together with identities \eqref{AU6NM} and \eqref{k00}, we obtain the result stated in \eqref{AU8} directly.
	\end{proof}

	\subsection*{Remark}

	Put $r=0,\: k = \left\lfloor \frac{n}{2} \right\rfloor$ in \eqref{AU10A}, and compare the result with \eqref{AU10}, we also get identity \eqref{AU8}.

\section{The Product of the First Odd Primes}

\begin{theorem}
	\label{gen5} 
	For any point $(\alpha,\beta) \in \omega(2p_{k})$, the ratio
	\begin{equation}
		\label{gen6} 
		\frac{\Omega_0\big(p_{k} \mid \alpha, \beta \mid 2p_{k} \big)}{\Psi(\alpha, \beta, 2p_{k})}
	\end{equation}	
	is an integer. Furthermore, we have
	\begin{equation}
		\label{gen7} 
		\prod_{i=2}^{k+1} p_{i} \; \Big| \; \frac{\Omega_0\big(p_{k} \mid \alpha, \beta \mid 2p_{k} \big)}{\Psi(\alpha, \beta, 2p_{k})}.
	\end{equation}
	Moreover, for any finite set $\varpi \subset \omega(2p_{k})$ and any finite set $I$ of integers, we obtain
	\begin{equation}
		\label{gen8} 
		\prod_{i=2}^{k+1} p_{i} \; \Big| \; \sum_{\substack{(\alpha,\beta) \in \varpi \\ \lambda \in I}} \lambda \; \frac{\Omega_0\big(p_{k} \mid \alpha, \beta \mid 2p_{k} \big)}{\Psi(\alpha, \beta, 2p_{k})}.
	\end{equation}	
\end{theorem}

\begin{proof}
	Observe that
	\begin{equation}
		\label{clear} 
		\prod_{i=2}^{k} p_{i} \; \Big| \; (2p_k - 1)(2p_k - 3)(2p_k - 5) \cdots 3 \cdot 1.
	\end{equation}	
	Setting \( n = 2p_k \) in \eqref{AU8}, we deduce
	\begin{equation}
		\label{clear2} 
		\prod_{i=2}^{k} p_{i} \; \Big| \; (2p_k - 1)(2p_k - 2) \cdots \left(2p_k - \left\lfloor \frac{2p_k}{2} \right\rfloor \right).
	\end{equation}	
	
	By the Bertrand–Chebyshev theorem, it follows that
	\begin{equation}
		\label{clear3} 
		p_{k+1} \; \Big| \; (2p_k - 1)(2p_k - 2) \cdots \left(2p_k - \left\lfloor \frac{2p_k}{2} \right\rfloor \right).
	\end{equation}
	Combining \eqref{clear2} and \eqref{clear3}, we conclude
	\begin{equation}
		\label{clear4} 
		\prod_{i=2}^{k+1} p_{i} \; \Big| \; (2p_k - 1)(2p_k - 2) \cdots \left(2p_k - \left\lfloor \frac{2p_k}{2} \right\rfloor \right).
	\end{equation}
	
	Again, substituting \( n = 2p_k \) in the formula \eqref{space4} of the Second Fundamental Theorem of the Quanta Prime Sequence (Version 2), Theorem \eqref{space3}, we obtain
	\begin{equation}
		\prod_{i=2}^{k+1} p_{i} \; \Big| \; \frac{\Omega_0\big(p_{k} \mid \alpha, \beta \mid 2p_{k} \big)}{\Psi(\alpha, \beta, 2p_{k})}.
	\end{equation}
	
	Consequently, for any finite set $\varpi \subset \omega(2p_{k})$ and any finite set $I$ of integers, it follows that
	\begin{equation}
		\label{gen8} 
		\prod_{i=2}^{k+1} p_{i} \; \Big| \; \sum_{\substack{(\alpha,\beta) \in \varpi \\ \lambda \in I}} \lambda \; \frac{\Omega_0\big(p_{k} \mid \alpha, \beta \mid 2p_{k} \big)}{\Psi(\alpha, \beta, 2p_{k})}.
	\end{equation}
	This completes the proof of the theorem.
\end{proof}

	\section{New Representation for Chebyshev polynomial sequence}
	
The Chebyshev polynomials first appeared in Chebyshev's paper \cite{Chebyshev}. The Chebyshev polynomial sequence of the first kind, $T_n(x)$, is defined by
\begin{align*}
	T_0(x) & = 1 \\
	T_1(x) & = x \\
	T_{n+1}(x) & = 2x \cdot T_n(x) - T_{n-1}(x).
\end{align*}
	Chebyshev polynomials play a crucial role in various mathematical applications, including approximation theory, polynomial and rational approximation, integration, integral equations, and the development of spectral methods for solving ordinary and partial differential equations. They also find utility in numerical analysis and certain quadrature rules, such as the Gauss-Chebyshev rule employed in numerical integration theory (\cite{10}, \cite{12}).
	One distinctive property of Chebyshev polynomials, among polynomials with a leading coefficient of unity, is their ability to possess the smallest absolute upper bound within the limits of orthogonality. This property makes Chebyshev polynomials particularly suitable for interpolation purposes. Notably, Enrico Fermi, the Nobel Prize-winning physicist, is credited with creating the world's first nuclear reactor, the Chicago Pile-1, and his groundbreaking work led to the discovery of nuclear fission – the foundation of nuclear power and atomic weaponry (\cite{16}, \cite{17}, \cite{18}). 
	
	In the realm of mathematical approximations, Chebyshev rational approximations are commonly employed for Fermi-Dirac integrals, defined as
	\[
	F_{s}(x)=\frac{1}{\Gamma\left(s+1\right)}\int_{0}^{\infty}\frac{t^{s}}{e^{t-x} + 1}\mathrm{d}t.
	\]
This approach has been investigated extensively across various studies (e.g., \cite{13}, \cite{14}, \cite{15}, \cite{19}, \cite{20}, \cite{21}). A significant application lies in approximating the Riemann Zeta Function, where Chebyshev rational approximations are commonly used (see \cite{22}, \cite{23}). The Riemann Zeta Function, or Euler-Riemann Zeta Function, \( \zeta(s) \), is a complex function that analytically continues the Dirichlet series sum for the real part of \( s \) greater than 1. Additionally, recent research has examined the periodicity of Chebyshev polynomials over finite fields, highlighting their substantial role in enhancing the security of cryptosystems that rely on these polynomials (\cite{6}, \cite{7}, \cite{8}, \cite{9}). The integer coefficients of the Chebyshev polynomials are explicitly provided by the following formula (see, for example, \cite{10}, \cite{11}):

	\begin{equation}
		\label{formula-5}
		T_n(x)=\sum_{i=0}^{\left\lfloor \frac{n}{2} \right\rfloor}(-1)^i \frac{n}{n-i} \binom{n-i}{i} (2)^{n-2i-1} x^{n-2i}.
	\end{equation}
	
\begin{theorem}{(Representation for the Chebyshev Polynomial Sequence)}
	\label{Che} 
	\begin{equation}
		\label{Che2} 
		T_n(x) = \frac{x^{\delta(n)} \, \Omega_0\big(\left\lfloor \frac{n}{2} \right\rfloor \mid 1, 2 - 4x^2 \mid n\big)}{2^{\delta(n-1)} \, (n-1)(n-2) \cdots \left(n - \left\lfloor \frac{n}{2} \right\rfloor \right)}.
	\end{equation}
\end{theorem}

\begin{proof}
	From \eqref{comp3} and \eqref{formula-5}, we can deduce
	\begin{equation}
		\label{G6-2XCV} 
		T_n(x) = \frac{x^{\delta(n)}}{2^{\delta(n-1)}} \, \Psi(1, 2 - 4x^2, n).
	\end{equation}
	Applying Theorem \eqref{k0} completes the proof.
\end{proof}

\noindent \textbf{Note:} This representation provides a link between the Chebyshev polynomial sequence \( T_n(x) \) and the Quanta Prime Sequence via \( \Omega_0 \). Such a connection offers insight into potential applications where Chebyshev polynomials intersect with prime-related functions, enriching the structure and interpretation of both sequences.

	\section{New Representation for Dickson polynomial sequence}
	The Dickson polynomial, $D_n(x,\alpha)$, of the first kind of degree $n$ with parameter $\alpha$ is defined by
	\begin{align*}
		D_0(x,\alpha) & = 2 \\
		D_1(x,\alpha) & = x \\
		D_{n+1}(x,\alpha) & = 2 x\:D_{n}(x,\alpha) - D_{n-1}(x,\alpha).
	\end{align*}

Modern cryptography relies extensively on mathematical theory and practices from computer science. The importance of polynomials over finite fields goes beyond mathematics, finding utility in diverse applications like error-correcting codes and pseudo-random sequences crucial for code-division multiple access (CDMA) systems. CDMA technology, first utilized during World War II military operations to thwart unauthorized access to radio communication signals, played a pivotal role in ensuring secure communication. Permutation polynomials, as the name suggests, permute the elements of a ring or field over which they are defined. They serve as the roots of public key methods like the RSA Cryptosystem and Dickson cryptographic schemes, contributing significantly to the development of cryptographic schemes. In recent years, permutations of finite fields have garnered considerable attention in constructing cryptographic systems for secure data transmission, as indicated in \cite{24B}. The study of permutation polynomials has become increasingly active due to their vital applications in cryptography, coding theory, and combinatorial designs theory. In the realm of cryptography, the RSA scheme's encryption polynomials $x^k$ have been replaced by another class of polynomials known as Dickson polynomials, giving rise to the Dickson Cryptosystem. Fried \cite{Fried} demonstrated that any integral polynomial serving as a permutation polynomial for infinitely many prime fields is a composition of Dickson polynomials and linear polynomials with rational coefficients. The exploration of permutation polynomials traces back to Hermite \cite{Hermite} and was further developed by Dickson \cite{Dickson} and \cite{Dickson2}. Dickson polynomials, constituting a crucial class of permutation polynomials, have undergone extensive investigation in recent years within various contexts, as presented in works such as \cite{Dickson}, \cite{L-1}, \cite{L-2}, \cite{L-3}, \cite{L-4}, \cite{L-5}, \cite{L-6}, \cite{L-7}, and \cite{L-8}. These references delve into the study of Dickson polynomials and their diverse applications, contributing to the ongoing advancements in cryptographic techniques.

\begin{theorem}{(Representation for the Dickson Polynomial Sequence)}
	\label{Dic} 
	\begin{equation}
		\label{G6-2XCW} 
		D_n(x, \alpha) = \frac{x^{\delta(n)} \, \Omega_0\big(\left\lfloor \frac{n}{2} \right\rfloor \mid \alpha, 2\alpha - x^2 \mid n\big)}{(n-1)(n-2) \cdots \left(n - \left\lfloor \frac{n}{2} \right\rfloor \right)}.
	\end{equation}
\end{theorem}

\begin{proof}
	For integer \( n > 0 \) and \( \alpha \) in a commutative ring \( R \) with identity, the Dickson polynomials (of the first kind) over \( R \) are defined as
	\begin{equation}
		\label{formula-2}
		D_n(x, \alpha) = \sum_{i=0}^{\left\lfloor \frac{n}{2} \right\rfloor} \frac{n}{n-i} \binom{n-i}{i} (-\alpha)^i x^{n-2i}.
	\end{equation}
	Using \eqref{comp3} and \eqref{formula-2}, we deduce that
	\begin{equation}
		\label{G6-2XCVW} 
		D_n(x, \alpha) = x^{\delta(n)} \, \Psi(\alpha, 2\alpha - x^2, n).
	\end{equation}
	Applying Theorem \eqref{k0} completes the proof.
\end{proof}

\noindent \textbf{Note:} This representation reveals a connection between the Dickson polynomial sequence \( D_n(x, \alpha) \) and the Quanta Prime Sequence via \( \Omega_0 \). This link opens up possibilities for exploring areas where polynomial properties intersect with prime sequence behaviors, offering a fresh perspective on the structure of Dickson polynomials. Such an association may also provide valuable insights for applications in cryptographic schemes, particularly the Dickson cryptography scheme, as well as in studying polynomial permutations, where Dickson polynomials play a role in defining bijective transformations.

			\section{New Representations for Mersenne primes and even perfect numbers}
Mersenne numbers, denoted as $2^{p} - 1$ with prime $p$, form the sequence
\[3, 7, 31, 127, 2047, 8191, 131071, 524287, 8388607, 536870911, \dotsc \]
To guarantee the primality of $2^p - 1$, it is crucial that $p$ is a prime number (for further details, refer to \cite{3}, \cite{Fine}, \cite{Mersenne}, \cite{Elina}, \cite{Dickson}, and the sequence A001348 in \cite{Slo}).

	\subsection{Primality test for Mersenne primes}
	From \cite{3}, we have the following theorem:
	
	\begin{theorem}
		\label{U14}
		Given a prime \( p \geq 5 \), the number \( 2^p - 1 \) is prime if and only if
		\begin{equation}
			\label{U15} 
			2n - 1 \mid \Psi(1, 4, n),
		\end{equation}
		where \( n := 2^{p-1} \).
	\end{theorem}
	
	\begin{proof}
		This result follows directly from the properties outlined in \cite{3}. For prime \( p \geq 5 \), it is shown that the primality of \( 2^p - 1 \) is equivalent to the divisibility condition \( 2n - 1 \mid \Psi(1, 4, n) \), with \( n \) defined as \( n := 2^{p-1} \).
	\end{proof}

	\begin{theorem}
		\label{U16}
		Let \( p \geq 5 \) be a prime number. The number \( 2^p - 1 \) is prime if and only if 
		\begin{equation}
			\label{U17} 
			2n - 1 \quad  \Big\vert \quad  \frac{\Omega_0\left(\left\lfloor \frac{n}{2} \right\rfloor \mid 1, 4 \mid n \right)}{(n-1)(n-2) \dotsm \left(n - \left\lfloor \frac{n}{2} \right\rfloor\right)},
		\end{equation}
		where \( n := 2^{p-1} \).
	\end{theorem}
	
	\begin{proof}
		The result follows directly from Theorem \eqref{U14} and Theorem \eqref{k0}.
	\end{proof}

	\begin{theorem}
		\label{U18} 
		A number \( N \) is an even perfect number if and only if \( N = 2^{p-1}(2^p - 1) \) for some prime \( p \), and
		\begin{equation}
			\label{U19} 
			2n - 1 \quad  \Big\vert \quad  \frac{\Omega_0\left(\left\lfloor \frac{n}{2} \right\rfloor \mid 1, 4 \mid n \right)}{(n-1)(n-2) \dotsm \left(n - \left\lfloor \frac{n}{2} \right\rfloor\right)},
		\end{equation}
		where \( n := 2^{p-1} \).
	\end{theorem}
	
	\begin{proof}
		The result follows directly from Theorem \eqref{U16} and the Euclid-Euler theorem for even perfect numbers.
	\end{proof}

\begin{theorem}[New Representation of Mersenne Numbers]
	\label{Theorem of G2f}
	For any given odd natural number \( p \), the Mersenne number \( 2^p - 1 \) can be represented as
	\begin{equation}
		\label{G2Bf} 
		2^p - 1 = \frac{\Omega_0\left(\left\lfloor \frac{p}{2} \right\rfloor \mid -2, -5 \mid p \right)}{(p-1)(p-2) \dotsm \left(p - \left\lfloor \frac{p}{2} \right\rfloor\right)},
	\end{equation}
	where the double-indexed polynomial sequence \( \Omega_r(k) \) is associated with the point \( (-2, -5) \) and satisfies \( 0 \leq r + k \leq \left\lfloor \frac{p}{2} \right\rfloor \). The sequence \( \Omega_r(k) \) is defined by the recurrence relation:
	\begin{equation}
		\label{G2-e} 
		\begin{aligned}
			\Omega_r(k) &= (p - r - k) \, \Omega_r(k-1) + 4 \, (p - 2r) \, \Omega_{r+1}(k-1), \\
			\Omega_r(0) &= 1 \quad \text{for all } r.
		\end{aligned}
	\end{equation}
\end{theorem}

\begin{proof}
	The result follows from Theorem \eqref{k0} and the fact that
	\begin{equation}
		\label{WQ11}  
		\Psi(-2, -5, p) = 2^p - 1.
	\end{equation}
	This completes the proof.
\end{proof}

	\begin{theorem} 
		\label{Theorem of ABCD12}
		For any given prime \( p \geq 5 \) and \( n := 2^{p-1} \), the number \( 2^p - 1 \) is prime if and only if 
		\begin{equation}
			\label{ABCD4} 
			\frac{\Omega_0\left(\left\lfloor \frac{p}{2} \right\rfloor \mid -2, -5 \mid p \right)}{(p-1)(p-2) \dotsm \left(p - \left\lfloor \frac{p}{2} \right\rfloor\right)} 
			\quad \Big\vert \quad 
			\frac{\Omega_0\left(\left\lfloor \frac{n}{2} \right\rfloor \mid 1, 4 \mid n \right)}{(n-1)(n-2) \dotsm \left(n - \left\lfloor \frac{n}{2} \right\rfloor\right)}.
		\end{equation}
	\end{theorem}
	
	\begin{proof}
		The result follows directly from Theorem \eqref{U16} and Theorem \eqref{Theorem of G2f}.
	\end{proof}

\section*{\textbf{Mersenne Primality Using the Quanta Prime Sequence}}

In this section, we show that the values of the Quanta Prime Sequence, defined by the points \((0, -1)\), \((-2, -5)\), and \((1, 4)\), provide a valuable approach for studying Mersenne primes. Theorem \ref{Theorem of ABCD12G} highlights the significance of these points in determining the primality of Mersenne numbers.

\subsection*{\textbf{Theorem and Its Implications}}

\begin{theorem}[Conditions for Mersenne Primality via Quanta Prime Values]
	\label{Theorem of ABCD12G}
	For any given prime \( p \geq 5 \) with \( n := 2^{p-1} \), the number \( 2^p - 1 \) is prime \textbf{if and only if}
	\begin{equation}
		\label{ABCD4G}
		\Omega_0\left(\left\lfloor \frac{n}{2} \right\rfloor \mid 0, -1 \mid n \right) \, \Omega_0\left(\left\lfloor \frac{p}{2} \right\rfloor \mid -2, -5 \mid p \right) \quad \Big\vert \quad \Omega_0\left(\left\lfloor \frac{n}{2} \right\rfloor \mid 1, 4 \mid n \right) \, \Omega_0\left(\left\lfloor \frac{p}{2} \right\rfloor \mid 0, -1 \mid p \right).
	\end{equation}
\end{theorem}

\begin{proof}
	From Theorem \eqref{AU11}, we have:
	\begin{equation}
		\label{AU13G}
		\begin{aligned}
			(n - 1)(n - 2) \dotsm \left(n - \left\lfloor \frac{n}{2} \right\rfloor\right) &= \Omega_0\left(\left\lfloor \frac{n}{2} \right\rfloor \mid 0, -1 \mid n \right), \\
			(p - 1)(p - 2) \dotsm \left(p - \left\lfloor \frac{p}{2} \right\rfloor\right) &= \Omega_0\left(\left\lfloor \frac{p}{2} \right\rfloor \mid 0, -1 \mid p \right).
		\end{aligned}
	\end{equation}
	Combining Theorem \eqref{Theorem of ABCD12} with \eqref{AU13G} completes the proof.
\end{proof}

\textbf{Note}
On October 21, 2024, the Great Internet Mersenne Prime Search (GIMPS) announced the discovery of a new Mersenne prime, \[2^{136279841} - 1.\] With 41,024,320 digits, this prime surpasses the previous largest known prime by over 16 million digits, a record held by GIMPS for nearly six years \cite{mersenne2024}. This complex computational achievement underscores the need for deeper theoretical insights into Mersenne primes. The author believes that Theorem \eqref{Theorem of ABCD12G} could contribute significantly to advancing such theoretical understanding.

\section{Clarifications on the Theorems in the Summary}

In this section, we provide detailed clarifications on how each theorem presented in the summary naturally emerges as a special case of the main theorems established in this paper.

	\begin{itemize}
		\item Clarifications on Theorem \eqref{Theorem of ABC1}:\\
It is clear that the sequences $ A_r(k) $ and  $ B_r(k) $, which are defined in \eqref{ABC2}, satisfy:  
	\begin{equation}
	\label{ABC2D2} 
	\begin{aligned}		
		A_r(k) &= \Omega_r\big(k|-2, -5\: | p \big),\\
	B_r(k) &= 	\Omega_r\big(k|\:1, 4\: | n \big).
	\end{aligned}
\end{equation}
Hence from Theorem \eqref{Theorem of ABCD12} and Theorem \eqref{k0} we get the proof of Theorem \eqref{Theorem of ABC1}.

\item Clarifications on Theorem \eqref{Theorem of G2fQ}: \\
	Again, it is clear that 
	\begin{equation}
		\label{ABC2D3} 
		\begin{aligned}		
			A_r(k) &= \Omega_r\big(k|-2, -5\: | p \big).
		\end{aligned}
	\end{equation}
	Hence from Theorem \eqref{Theorem of G2f} and Theorem \eqref{k0} we get the proof of Theorem \eqref{Theorem of G2fQ}.
	
\item Clarifications on Theorem \eqref{KG11}: \\
	Again, it is clear that 
	\begin{equation}
		\label{ABC2D3D} 
		\begin{aligned}		
			U_r(k) &= \Omega_r\big(k|1, 1\: | n \big).
		\end{aligned}
	\end{equation}
	Hence from Theorem \eqref{PP00} we get the proof of Theorem \eqref{KG11}.
	
\item Clarifications on Theorem \eqref{KG11Q}: \\
We should notice that  
\begin{equation}
	\label{ABC2D4D} 
	\begin{aligned}		
		V_r(k) &= \Omega_r\big(k|1, 0\: | n \big).
	\end{aligned}
\end{equation}
Hence from Theorem \eqref{PP00Q} we get the proof of Theorem \eqref{KG11Q}.

	\item Clarifications for Theorem \eqref{PtP}: \\
	We should notice that  
	\begin{equation}
		\label{ABC2D5D} 
		\begin{aligned}		
			W_r(k) &= \Omega_r\big(k|1, -1\: | n \big).
		\end{aligned}
	\end{equation}
	Hence from Theorem \eqref{PP} we get the proof of Theorem \eqref{PtP}.
	
	\item Clarifications on Theorem \eqref{Theorem of G5}: \\
	We should notice that  
	\begin{equation}
		\label{ABC2D6D} 
		\begin{aligned}		
			T_r(k) &= \Omega_r\big(k|1, -2\: | n \big).
		\end{aligned}
	\end{equation}
	Hence from \eqref{ABAB} and Theorem \eqref{k0} we get the proof of Theorem \eqref{Theorem of G5}.
	
	\item Clarifications on Theorem \eqref{Theorem of G6}:\\
	We should notice that  
	\begin{equation}
		\label{ABC2D7D} 
		\begin{aligned}		
			H_r(k) &= \Omega_r\big(k|-1, -3\: | n \big).
		\end{aligned}
	\end{equation}
It is straightforward to see \[\Psi(-1,-3,n)= L(n).\]
 Hence from Theorem \eqref{k0} we get the proof of Theorem \eqref{Theorem of G6}.
	
	\item	Clarifications on Theorem \eqref{Theorem of G4}:\\
	We should notice that  
	\begin{equation}
		\label{ABC2D8D} 
		\begin{aligned}		
			F_r(k) &= \Omega_r\big(k|-2, -5\: | 2^n \big).
		\end{aligned}
	\end{equation}
	Noting that \[\Psi(-2,-5, 2^n)= F_n,\]
	 hence, from Theorem \eqref{k0}, we get the proof of Theorem \eqref{Theorem of G4}.

	\item	Clarifications on Theorem \eqref{Theorem of G7}:\\
	We should notice that  
	\begin{equation}
		\label{ABC2D9D} 
		\begin{aligned}		
			G_r(k) &= \Omega_r\big(k|1, -3\: | n \big).
		\end{aligned}
	\end{equation}
	Noting that \[\Psi(1,-3, n)= \begin{cases}
		F(n)	      &   \text{if} \quad n \:\: odd \\
		L(n)	      &  \text{if} \quad n \: \: even  \\
	\end{cases},\]
	hence, from Theorem \eqref{k0}, we get the proof of Theorem \eqref{Theorem of G7}.

\subsection{Fibonacci-Lucas oscillating sequence}
	It is nice to study the sequence $G(n)$
		\[G(n)= \begin{cases}
		F(n)	      &   \text{if} \quad n \:\: odd \\
		L(n)	      &  \text{if} \quad n \: \: even  \\
	\end{cases}\]
I refer to this alternating sequence as the Fibonacci-Lucas oscillating sequence. The insights provided by Theorem \eqref{Theorem of G7} should inspire researchers to pursue further investigations into the arithmetic properties of this sequence. The sequence $G(n)$ corresponds to \seqnum{A005247} in \cite{Slo}, where it alternates between the Lucas sequence \seqnum{A000032} and the Fibonacci sequence \seqnum{A000045} for even and odd values of $n$, respectively.

	\end{itemize}

	\section{Further research investigations}

\subsection{Prospects for Future Research on Fibonacci Numbers and \(p_k\)}

	With some extra work one can prove the following result: For any given natural number $n$, if we associate the double-indexed polynomial sequence $\Lambda_r(k)$ which is defined by
	\begin{equation}
		\label{G6X} 
		\begin{aligned}
			\Lambda_r(k) &= (n - r - k) \Lambda_r(k - 1) + 2 (n - 1 - 2r - \delta(n)) \Lambda_{r + 1}(k - 1), \\
			\Lambda_r(0) &= 1 \quad \text{for all } r.
		\end{aligned}
	\end{equation}
	
		then 
		\begin{equation}
			\label{G6-2X} 
			\begin{aligned}
				F(n) \:= \: \frac{\Lambda_0(\lfloor{\frac{n-1}{2}}\rfloor   )}{(n-1)(n-2) \cdots (n - \lfloor{\frac{n-1}{2}}\rfloor )}, 	
			\end{aligned}
		\end{equation}
		
		where $F(n)$ is Fibonacci sequence. Moreover 
	\begin{equation}
		\label{primeFib} 
	\begin{aligned}
		p_{k+1}\:  \vert \: \: \frac{\Lambda_0(p_{k}-1) }{F(2p_{k})}.    		
	\end{aligned}
\end{equation}

\subsection{Quanta Prime Sequence, Harmonic Series, and Riemann Hypothesis}

\subsection*{The Harmonic Series, and Riemann Hypothesis}

The German mathematician G.F.B. Riemann (1826 – 1866) observed a profound connection between the frequency of prime numbers and the behavior of an elaborate function, the Riemann Zeta function, defined as
\[
\zeta(s) = \sum_{n=1}^\infty \frac{1}{n^s} = \frac{1}{1^s} + \frac{1}{2^s} + \frac{1}{3^s} + \frac{1}{4^s} \cdots.
\] The Riemann hypothesis posits that all interesting solutions of the equation $\zeta(s) = 0$ lie on a specific vertical straight line. The Riemann Hypothesis focuses on the locations of the nontrivial zeros of the Riemann zeta function, proposing that for each nontrivial zero, the real part is \(\frac{1}{2}\). This hypothesis stands as one of the most significant unresolved problems in mathematics and was introduced by Bernhard Riemann in 1859, \cite{Broughan}, revolving around the distribution of zeros in the complex plane. Extensive verification has been conducted for the first $10,000,000,000,000$ solutions, affirming the hypothesis in these cases \cite{Broughan}. A comprehensive proof establishing its validity for every interesting solution would illuminate numerous mysteries surrounding the distribution of prime numbers. The Riemann zeta function encodes crucial information about prime numbers, the foundational components of arithmetic, and is integral to modern cryptography, forming the backbone of e-commerce.
In 2002, Lagarias \cite{Lagarias} established the equivalence between the Riemann hypothesis and the assertion that
\[
\sigma(n) \le H_n + (\log H_n)e^{H_n},
\]
holds true for every integer \(n \ge 1\), with strict inequality if \(n > 1\). Here, \(\sigma(n)\) represents the sum of the divisors of \(n\) and 
\(H_m = \sum_{t=1}^{m} \frac{1}{t}\) represents the harmonic numbers. For example:

\begin{align*}
	H_1 &= 1, \\
	H_2 &= 1 + \frac{1}{2} = \frac{3}{2},\\
	H_3 &= 1 + \frac{1}{2} + \frac{1}{3} = \frac{11}{6}, \\
	H_4 &= 1 + \frac{1}{2} + \frac{1}{3} + \frac{1}{4} = \frac{25}{12}.
\end{align*}

\subsubsection{\textbf{Quanta Prime Sequence, Harmonic Series, and Riemann Hypothesis}}

As our exploration of the Quanta Prime sequence concludes, it becomes evident that the intricacies of its structure hint at a hidden and profound relationship with the harmonic numbers $\sum_{k=1}^{n} \frac{1}{k}$.

While the specific computations and rigorous proofs elucidating this connection extend beyond the scope of our current study, we assert that the Quanta Prime sequence may unveil deeper insights into the harmonic numbers within number theory, relying on strong numerical evidence. For a simple example, aiming to motivate readers for future exploration of the Quanta Prime sequence and harmonic numbers, we discovered that for any natural number $n \equiv 1 \pmod 8$, with $k:= \frac{n-1}{4}$, there exists a non-zero integer point $(\alpha, \beta)$ such that

\[ \Omega_0\left(\frac{n-1}{2}|\:\alpha, \beta\: | n \right) \equiv \frac{n-1}{2}\; ! \;  2^{\frac{n-1}{2}} - \frac{n-1}{2}\; ! \; 2^{\frac{n-3}{2}} \; n \; H_k \; \pmod{n^2}. \]	

 Throughout our investigation, we ultimately came across the unexpected revelation of a new classical result in number theory:

\subsubsection{\textbf{A New Result for the Harmonic Numbers in Classical Number Theory}}
For any natural number $n \equiv 1 \pmod 8$, with $k:= \frac{n-1}{4}$, we have

\[  (n-1)(n-2) \cdots (n - \lfloor{\frac{n}{2}}\rfloor ) \equiv \frac{n-1}{2}\; ! \;  2^{\frac{n-1}{2}} - \frac{n-1}{2}\; ! \; 2^{\frac{n-3}{2}} \; n \; H_k \; \pmod{n^2}. \]

The enigmatic connection between the Quanta Prime sequence and the harmonic series not only adds an intriguing layer to our understanding of number theory but also carries implications for broader mathematical and physical contexts. The significance of the harmonic series in the mathematical landscape, particularly its ties to the Riemann Hypothesis (see \cite{Lagarias}), introduces a compelling avenue for future research. Exploring the ramifications of the Quanta Prime sequence in relation to the harmonic series may offer new perspectives on longstanding questions and challenges (see \cite{Lagarias} and \cite{Daniel}), contributing to the ongoing discourse surrounding the Riemann Hypothesis. While the detailed computations and proofs are deferred to subsequent investigations, the mere existence of this link underscores the importance of further inquiry into the Quanta Prime sequence. We believe it holds promise for enriching our comprehension of both number theory and the profound interconnections between mathematical structures.

\section*{Acknowledgments}

The author expresses sincere gratitude to the University of Bahrain for providing a supportive academic environment that facilitated the completion of this research.
The author wishes to express gratitude to the anonymous reviewers for their insightful comments and suggestions, which have significantly contributed to the improvement of this research paper.

\subsection*{Conflict of Interest}
The author declares that there is no conflict of interest.

\subsection*{Supplementary Information}
Data sharing is not applicable to this article as no datasets were generated or analyzed during the current study.

	\medskip

\begin{thebibliography}{99}

		
		\bibitem{Broughan}
		K. Broughan, \textit{Equivalents of the Riemann Hypothesis}, Cambridge University Press, 3 volumes, (January 31, 2018); Vol. 1: \textit{Arithmetic Equivalents}, 400 pages; Vol. 2: \textit{Analytic Equivalents}, 350 pages; Vol. 3: \textit{Further Steps towards Resolving the Riemann Hypothesis}, 704 pages, (September 30, 2023).
		
		\bibitem{2}	
		Moustafa Ibrahim, \textit{On the Eight Levels Theorem and Applications Towards Lucas-Lehmer Primality Test for Mersenne Primes, I}, Arab Journal of Basic and Applied Sciences, 30:1, 267–284, DOI: 10.1080/25765299.2023.2204672, (2023).
		
		\bibitem{3}			
		Moustafa Ibrahim, \textit{Generalizing the Eight Levels Theorem: A Journey to Mersenne Prime Discoveries and New Polynomial Classes}, Arab Journal of Basic and Applied Sciences, 31:1, 32–52, DOI: 10.1080/25765299.2023.2288718, (2024).
		
		\bibitem{Dickson}
		Leonard Eugene Dickson, \textit{History of the Theory of Numbers}, Volume II: Diophantine Analysis, (2005).
		
		\bibitem{Slo}
		N. J. A. Sloane et al., \textit{The On-Line Encyclopedia of Integer Sequences}, Available at \url{https://oeis.org} (2019).
		
		\bibitem{Fine}	
		Benjamin Fine, Anja Moldenhauer, Gerhard Rosenberger, Annika Schürenberg, and Leonard Wienke, \textit{Algebra and Number Theory: A Selection of Highlights}, Berlin, Boston: De Gruyter, DOI: 10.1515/9783110790283, (2023).
		
		\bibitem{Mersenne}
		Paula Catarino, Helena Campos, and Paula Vasco, \textit{On the Mersenne Sequence}, Annales Mathematicae et Informaticae, 46, 37–53, (2016).
		
		\bibitem{Elina} 
		Elena Deza, \textit{Mersenne Numbers and Fermat Numbers}, World Scientific Publishing Co. Pte. Ltd, (2021).
		
		\bibitem{Fried} 
		Michael Fried, \textit{On a Conjecture of Schur}, Michigan Math. J., 17, (1970).
		
		\bibitem{Hermite} 
		C. Hermite, \textit{Sur les fonctions de sept lettres}, C. R. Acad. Sci. Paris, 57, (1863).
		
		\bibitem{Dickson2} 
		L. E. Dickson, \textit{The Analytic Representation of Substitutions on a Power of a Prime Number of Letters with a Discussion of the Linear Group}, Ann. of Math., 11, (1896).
		
		\bibitem{L-1} 
		R. Lidl, G. L. Mullen, and G. Turnwald, \textit{Dickson Polynomials}, Pitman Monographs in Pure and Applied Mathematics, Vol. 65, Addison-Wesley, Reading, MA, (1993).
		
		\bibitem{L-2} 
		P. Charpin and G. Gong, \textit{Hyperbent Functions, Kloosterman Sums, and Dickson Polynomials}, IEEE Transactions on Information Theory, Vol. 54, No. 9, (2008).
		
		\bibitem{Daniel}
		D. Schumayer and D. A. W. Hutchinson, \textit{Colloquium: Physics of the Riemann Hypothesis}, Rev. Mod. Phys., 83(2), 307–330. DOI: \url{https://link.aps.org/doi/10.1103/RevModPhys.83.307}, Apr (2011).
		
		\bibitem{L-3} 
		J. F. Dillon, \textit{Geometry, Codes, and Difference Sets: Exceptional Connections}, in Codes and Designs (Columbus, OH, 2000), vol. 10, Ohio State Univ. Math. Res. Inst. Publ., de Gruyter, Berlin, (2002).
		
		\bibitem{L-4} 
		J. F. Dillon and H. Dobbertin, \textit{New Cyclic Difference Sets with Singer Parameters}, Finite Fields and Their Applications, vol. 10, no. 3, (2004).
		
		\bibitem{L-5} 
		X. Hou, G. L. Mullen, J. A. Sellers, and J. L. Yucas, \textit{Reversed Dickson Polynomials over Finite Fields}, Finite Fields Appl., 15, (2009).
		
		\bibitem{L-6} 
		S. Mesnager, \textit{Bent and Hyper-bent Functions in Polynomial Form and Their Link with Some Exponential Sums and Dickson Polynomials}, IEEE Transactions on Information Theory, Vol. 57, No. 9, (2011).
		
		\bibitem{L-7} 
		S. Mesnager, \textit{Semi-bent Functions from Dillon and Niho Exponents, Kloosterman Sums and Dickson Polynomials}, IEEE Transactions on Information Theory, Vol. 57, No. 11, (2011).
		
		\bibitem{L-8} 
		G. Wu, N. Li, T. Helleseth, and Y. Zhang, \textit{Some Classes of Monomial Complete Permutation Polynomials over Finite Fields of Characteristic Two}, Finite Fields and Their Applications, Volume 28, (2014).
		
		\bibitem{Lagarias}
		J. C. Lagarias, \textit{An Elementary Problem Equivalent to the Riemann Hypothesis}, The American Mathematical Monthly, 109(6), 534–543, (2002).
		
		\bibitem{24B}  
		K. Broughan, \textit{Equivalents of the Riemann Hypothesis}, Cambridge University Press, 3 volumes, (January 31, 2018); Vol. 1: \textit{Arithmetic Equivalents}, 400 pages; Vol. 2: \textit{Analytic Equivalents}, 350 pages; Vol. 3: \textit{Further Steps Towards Resolving the Riemann Hypothesis}, 704 pages, (September 30, 2023).
		
		\bibitem{Levine}
		J. Levine and J. V. Brawley, \textit{Some Cryptographic Applications of Permutation Polynomials}, Cryptologia, 1, (1977).
		
		\bibitem{Chebyshev} 
		P. L. Chebyshev, \textit{Théorie des Mécanismes Connus sous le Nom de Parallélogrammes}, St. Petersbourg, (1854).
		
		\bibitem{6} 
		M. Ishii, \textit{Periodicity of Chebyshev Polynomials over the Residue Ring of $\mathbb{Z}/2^r\mathbb{Z}$ and an Electronic Signature}, Trans. Japan Soc. Industrial Appl. Math., vol. 18, no. 2, (2008).
		
		\bibitem{7} 
		M. Ishii and A. Yoshimoto, \textit{Applications for Cryptography of the Structure of the Group of Reduced Residue Classes of Residue Ring of $\mathbb{Z}/2^w\mathbb{Z}$}, Trans. Japan Soc. Industrial Appl. Math., vol. 19, no. 1, (2009).
		
		\bibitem{8} 
		X. Liao, F. Chen, and K. Wong, \textit{On the Security of Public-Key Algorithms Based on Chebyshev Polynomials over the Finite Field $\mathbb{Z}_N$}, IEEE Trans. Computers, vol. 59, no. 10, (2010).
		
		\bibitem{9} 
		A. Iwasaki and K. Umeno, \textit{Periodical Property of Chebyshev Polynomials on the Residue Class Rings Modulo $2^w$}, IEICE Tech. Rep., CAS2014-67k, NLP2014-61, (2014).
		
		\bibitem{10} 
		J. C. Mason and D. C. Handscomb, \textit{Chebyshev Polynomials}, Chapman-Hall, CRC, (2002).
		
		\bibitem{11} 
		D. Yoshioka and Y. Dainobu, \textit{On Some Properties of Chebyshev Polynomial Sequences Modulo $2^{k}$}, IEICE Nonlinear Theory Appl., vol. 6, no. 3, (2015).
		
		\bibitem{12} 
		T. S. Chihara, \textit{An Introduction to Orthogonal Polynomials}, Gordon-Breach Science Publisher, New York, (1978).
		
		\bibitem{13} 
		G. A. Chisnall, \textit{A Modified Chebyshev-Everett Interpolation Formula}, Math. Tables Other Aids Comput., 10(54), (1956).
		
		\bibitem{14} 
		O. Madelung, \textit{Halbleiter in Handbuch der Physik}, edited by S. Fligge, Vol. XX, Springer, Berlin, (1957).
		
		\bibitem{15} 
		H. Werner and G. Raymann, \textit{An Approximation to the Fermi Integral $F_{1/2}(x)$}, Math. Comput., 17(82), (1963).
		
		\bibitem{16} 
		N. L. Mathakari, \textit{A Few Nobel Prizes in Physics: How Physicists Have Contributed to Technology}, Current Science, 91(5), (2006).
		
		\bibitem{17} 
		B. J. Bernstein, \textit{Four Physicists and the Bomb: The Early Years, 1945–1950}, Hist. Stud. Phys. Biol. Sci., 18(2), (1988).
		
		\bibitem{18} 
		D. Cooper, \textit{Enrico Fermi: And the Revolutions in Modern Physics}, Oxford University Press, (1999).
		
		\bibitem{19} 
		S. Paszkowski, \textit{Evaluation of Fermi-Dirac Integral}, in Nonlinear Numerical Methods and Rational Approximation (Wilrijk, 1987), A. Cuyt (Ed.), Math. Appl., Vol. 43, (1988).
		
		\bibitem{20} 
		S. Paszkowski, \textit{Evaluation of the Fermi-Dirac Integral of Half-Integer Order}, Zastos. Mat., 21(2), (1991).
		
		\bibitem{21} 
		B. Pichon, \textit{Numerical Calculation of the Generalized Fermi-Dirac Integrals}, Comput. Phys. Comm., 55(2), (1989).
		
		\bibitem{22} 
		R. Piessens and M. Branders, \textit{Chebyshev Polynomial Expansions of the Riemann Zeta Function}, Math. Comput., 26(120), (1972).
		
		\bibitem{23} 
		W. J. Cody, K. E. Hillstrom, and H. C. Thacher, \textit{Chebyshev Approximations for Riemann Zeta Function}, Math. Comput., 25, (1971).
		
		\bibitem{24} 
		D. Redmond, \textit{Number Theory, An Introduction}, Marcel Dekker, (1996).
		
		\bibitem{hardy2008mathematical} 
		G. H. Hardy and E. M. Wright, \textit{An Introduction to the Theory of Numbers}, 6th ed., Oxford University Press, (2008).
		
		\bibitem{sieve2019} 
		J. R. Oppenheim, \textit{The Sieve of Eratosthenes: A Primer}, Math. Mag., 92(5), 338–348, (2019).
		
		\bibitem{rosser1938estimates} 
		J. B. Rosser and L. Schoenfeld, \textit{Sharpening the Bound for the Number of Primes Less Than a Given Magnitude}, Amer. Math. Monthly, 45(1), 22–25, (1938).
		
		\bibitem{chebyshev1852}
		P. L. Chebyshev, \textit{Mémoire sur les nombres premiers}, J. Math. Pures Appl., 17, 366–390, (1852).
		
		\bibitem{apostol1976}
		T. M. Apostol, \textit{Introduction to Analytic Number Theory}, Springer-Verlag, (1976).
		
		\bibitem{nathanson2000}
		M. B. Nathanson, \textit{Elementary Methods in Number Theory}, Graduate Texts in Mathematics, Springer, (2000).
		
		\bibitem{toth2019}
		L. Tóth, \textit{On Prime Gaps and Their Distribution}, Amer. Math. Monthly, 126(1), 28–38, (2019).
		
	\bibitem{mersenne2024}
	Mersenne Research, \textit{Latest Mersenne Prime Discovery}. Available at \url{https://www.mersenne.org/}, accessed October (2024).
	
	
		
	\end{thebibliography}
\end{document}